\documentclass[11pt,letterpaper]{article}

\usepackage[dvipsnames,svgnames,x11names]{xcolor}
\usepackage[final]{changes}
\definechangesauthor[color=BrickRed]{GB}
\definechangesauthor[color=NavyBlue]{JP}
\usepackage{todonotes}
\newcommand{\note}[2][]{\added[#1,remark={#2}]{}}

\usepackage{titling}
\usepackage[flushleft,para]{threeparttable}
\usepackage[toc,page]{appendix}
\usepackage{hyperref}
\hypersetup{
    colorlinks=true,
    linkcolor=blue,
    filecolor=magenta,      
    urlcolor=cyan,
    citecolor= blue
}

\usepackage[onehalfspacing]{setspace}
\usepackage{url}
\usepackage[affil-it]{authblk}
\usepackage{float}
\usepackage{makeidx}
\usepackage{adjustbox}
\usepackage{bookmark}
\usepackage{yhmath}
\usepackage{multirow}
\usepackage{amsfonts}
\usepackage{amsmath}
\usepackage{amssymb}
\usepackage{amsthm}
\usepackage[round]{natbib}
\usepackage{booktabs}
\usepackage{enumerate}
\usepackage{enumitem}
\usepackage{bm}
\usepackage{mathtools}
\usepackage{relsize}
\usepackage{comment} 
\usepackage{mathrsfs}
\usepackage{amsbsy}
\usepackage[scr=boondoxo,scrscaled=1.05]{mathalfa}
\usepackage{xfrac}
\usepackage{lineno}
\usepackage{graphicx}
\usepackage{fullpage}
\usepackage{setspace}
\allowdisplaybreaks[4]
\allowbreak

\newcommand{\q}{\mathbf{q}}

\newcommand{\e}[1]{\mathbb{E} \left[ #1 \right]}
\newcommand{\ee}[2]{\mathbb{E}_{#1} \left[ #2 \right]}

\newcommand{\vartt}{\texttt} %

\newcommand{\pbf}{\mathbf{p}}

\newcommand{\cX}{\mathcal{X}}

\newtheorem{proposition}{Proposition}

\usepackage{algorithm,algorithmicx,algpseudocode}
\algnewcommand{\Initialize}[1]{%
  \State \textbf{initialization:}
}
\algnewcommand{\STEPZERO}[1]{%
  \State \textbf{\underline{STEP 0. Initialization:}}
  \Statex \hspace*{\algorithmicindent}\parbox[t]{.8\linewidth}{\raggedright #1}
}
\algnewcommand{\STEPFIRST}[1]{%
  \State \textbf{\underline{STEP 2. Forward pass:}}
}
\algnewcommand{\STEPSECOND}[1]{%
  \State \textbf{\underline{STEP 1. Update lower bound:}}
}
\algnewcommand{\STEPTHIRD}[1]{%
  \State \textbf{\underline{STEP 3. Update upper bound:}}
}
\algnewcommand{\STEPFOURTH}[1]{%
  \State \textbf{\underline{STEP 4. Check stopping criterion:}}
}
\algnewcommand{\STEPFIFTH}[1]{%
  \State \textbf{\underline{STEP 5. Backward pass:}}
}
\makeatletter   %
\let\orgdescriptionlabel\descriptionlabel  %
\renewcommand*{\descriptionlabel}[1]{%
\let\orglabel\label   %
\let\label\@gobble   %
\phantomsection   %
\edef\@currentlabel{#1}%
\let\label\orglabel   %
\orgdescriptionlabel{#1}%
} %
\makeatother  %

 \usepackage{bibunits}
\defaultbibliographystyle{apalike}

\begin{document}

\title{A Multistage Distributionally Robust Optimization Approach to Water Allocation under Climate Uncertainty}

\author[1]{Jangho Park\thanks{\texttt{JanghoPark@lbl.gov}}}
\affil[1]{{\small Computational Research Division, Lawrence Berkeley National Laboratory, 1 Cyclotron Rd, Berkeley, CA 94720, United States}}
\author[2]{G\"{u}zin Bayraksan\thanks{\texttt{bayraksan.1@osu.edu}; Corresponding author}}
\affil[2]{{\small Department of Integrated Systems Engineering, The Ohio State University, 1971 Neil Avenue, Columbus, OH 43210, United States}}

\date{}

\maketitle
\vspace*{-0.7in}
	\begin{abstract}		
This paper investigates a Multistage Distributionally Robust Optimization (MDRO) approach to water allocation under climate uncertainty.  The MDRO is formed by creating sets of conditional distributions (called {\it conditional ambiguity sets}) on a  finite scenario tree. The distributions in the conditional ambiguity sets remain close to a nominal conditional distribution according a $\phi$-divergence (e.g., Kullback-Liebler divergence, Hellinger distance, Burg entropy, etc.).  The paper discusses a decomposition algorithm to solve the resulting MDRO and applies the modeling and solution techniques to allocate water in a rapidly-developing area of Tucson, Arizona. Tucson, like many arid and semi-arid regions around the world, faces considerable uncertainty in its ability to provide water for its citizens in the future. The primary sources of uncertainty in the Tucson region include (1) unpredictable population growth, (2) the availability of water from the Colorado River, and (3) the effects of climate variability on water consumption. This paper integrates forecasts for all these sources of uncertainty into a single optimization model for robust and sustainable water allocation. Then, it uses this model to analyze the value of constructing additional treatment facilities to reduce future water shortages. The results indicate that the MDRO approach can be very valuable for water managers by providing insights to minimize their risks and help them plan for the future.\smallskip

	\end{abstract}

\newpage
\section{Introduction}
\label{sec: intro_mdro}
\vspace*{-0.1in}

This paper studies Multistage Distributionally Robust Optimization (MDRO) with $\phi$-divergences and proposes a decomposition algorithm to solve an equivalent formulation of MDRO. It then applies the MDRO modeling and solution techniques for sustainable water allocation in a developing area of Tucson, Arizona.

More than 60\% of the water in Tucson is provided by the Colorado River. 
Without this water source, citizens of Tucson---as well as millions in Arizona, California, Nevada, and Mexico---are threatened. 
The Colorado River has been facing extreme water shortages in recent years.  
In 2015 and 2016, Lake Mead water elevation hit back-to-back record low, 
and, in June 2016, it reached its lowest level of 1071.64 feet for the first time in its 80-year history \citep{usbor_web}.  
As the Colorado River runs dry \citep{colorado} and population increases, the imbalance between supply and demand widens. 
And as climate variability threatens the Colorado River availability, it is imperative to sustainably manage this water resource by taking into account the many complex uncertainties it faces.

This paper presents a novel MDRO model with $\phi$-divergences to allocate Colorado River water to an urban water system in Tucson, AZ. Distributionally Robust Optimization (DRO) with $\phi$-divergences in the static/two-stage case has been proposed by the seminal work of \citet{bental2011robust}; see also further investigations by \cite{bayraksan2015data} and \cite{jiang_guan_16}. 
DRO acknowledges that uncertainties---like the long-term and complex ones on climate, population, and the Colorado River basin's hydrology considered in this paper---are not fully known. 
Such uncertainties, i.e., those whose probability distributions are not fully known, are referred to as {\it ambiguous} uncertainties.  
There is, however, historical data, sophisticated simulations, and detailed forecasts available from research, government, and utility sources. 
So, it is possible to build approximate future scenarios with an approximate nominal distribution for these ambiguous uncertainties. 
DRO then considers all distributions that are sufficiently close to this nominal distribution and optimizes a worst-case expected objective, where the expectations are taken with respect to all the considered distributions. %
The appeal of DRO is that it is more realistic because it explicitly considers existing data and forecasts, while acknowledging that these forecasts may contain errors. 

There is an increasing literature on various DRO formulations and solution techniques, but most of this literature focuses on static, two-stage, or chance-constrained settings  \citep{rahimian2019distributionally}. There is relatively little work on multistage DRO. Many of these works investigate different ways of forming a set of distributions in the multistage setting. Distance-based approaches consider staying sufficiently close to a nominal distribution according to a nested Wasserstein distance \citep{Pflug2014, Analui2014}, modified $\chi^2$ distance \citep{philpott2018}, $L_\infty$-norm \citep{Jianqiu2017Study} 
or Wasserstein-$\infty$ distance \citep{Dimitris2018Data}. An alternative approach uses moment-based sets \citep{xin2013time, xin2015distributionally, Bertsimas2019Adaptive, Babaei2019Data}, where the moments of the distributions must satisfy certain properties. \cite{Shapiro2016Rectangular} studies general MDRO theory, and \cite{shapiro2018tutorial} provides an overview. 

This paper, to the best of our knowledge for the first time, studies a general class of MDRO constructed by forming sets of conditional distributions via $\phi$-divergences on a finite-state, discrete-time stochastic process. It shows that many properties---such as risk aversion, decomposition, and reformulations like second-order cone representation---of the static/two-stage DRO extend to the multistage setting. Although many real-world problems (including the one studied in this paper) can only be appropriately analyzed using multistage models, such models are notoriously difficult to solve in general. To alleviate this difficulty, the paper presents a nested Benders decomposition algorithm \citep{Birge1985} 
to solve the resulting MDRO formed via a general class of $\phi$-divergences. Unlike most of the literature \citep[e.g.,][]{philpott2018, Jianqiu2017Study}, it does not assume independence (or simple forms of dependence) of the stochastic process between stages. However, it assumes a moderately large number of realizations so that MDRO can be solved in a reasonable time without resorting to sampling. Finally, the paper applies MDRO modeling framework to analyze a real-world water allocation problem and uses it to evaluate infrastructure decisions in the area. To the best of our knowledge, this is the first application of MDRO to an important operations research problem in natural resources considering climate variability.

Climate is one of the most important sources of uncertainty for long-term sustainability of water resources. 
Extensive research analyze the sensitivity of mitigation plans to uncertainties in climate  \citep[e.g.,][]{Singh_etal_14,harou2010economic,OHara2008,ROBERT20181033,DURANENCALADA2017567}. 
None of the existing work, however, considers the ambiguities in future climate predictions like the MDRO model. 
We believe the MDRO modeling framework is particularly important for problems with complex multiperiod uncertainties---like those that incorporate climate predictions.  %

One of the unique features of our model is that it combines various sources of data to generate water demand and supply scenarios.
We incorporate bias-corrected and spatially-downscaled global circulation climate models (formed via different organizations around the world), greenhouse gas concentration paths (as adopted by the Intergovernmental Panel on Climate Change (IPCC)), population forecasts (developed by governing agencies in the area), water-use trends as well as hydrological simulations of the Colorado River (conducted by the U.S. Bureau of Reclamation).
We explain our scenario generation methodology in Section \ref{sec: water}.

To cope with the possible future water shortages, we consider constructing additional water infrastructures.
New infrastructures cost hundreds of millions of dollars, so they should be evaluated carefully. We use the MDRO model to do so.
For many arid and semi-arid areas---like the area studied in this paper---{\it reclaimed water} (treated wastewater) is the only remaining water source \citep{woods2012centralized,lan_bayraksan_lansey_16}.  
Therefore, we consider constructing two decentralized water treatment facilities to increase reclaimed water use in the area. 
The first option builds a wastewater treatment plant and reuses treated wastewater for {\it nonpotable} (water that is used for many purposes except drinking) needs. This option saves freshwater resources for {\it potable} (water that is safe to drink) demands. 
The second alternative considers Indirect Potable Reuse (IPR), which treats wastewater to a very high quality and blends with other high-quality, drinkable water sources.

In summary, this paper presents a first MDRO approach for sustainable water allocation in urban water systems. 
It applies this model to allocate Colorado River water through mid-century to a developing area of Tucson, incorporating various uncertainties on climate, population, water-use trends, and the Colorado River water availability. 
This MDRO model is then used to assess water reuse strategies by evaluating the value of constructing additional water treatment facilities. 
It is important to highlight that the presented MDRO modeling and solution techniques are not limited to the water allocation problem studied in this paper. In fact, they have the potential to model and analyze many important operations research problems (e.g., in energy, finance) with substantial and complex multiperiod uncertainties, whose distributions cannot be fully known.

The rest of the paper is organized as follows.
Section \ref{sec: phi_divergences_water} reviews $\phi$-divergences and describes the MDRO with $\phi$-divergences.
This section also discusses a risk-averse interpretation of MDRO.
Section \ref{sec: alg} presents a decomposition algorithm to efficiently solve the MDRO.
The modeling and solution techniques of Sections \ref{sec: phi_divergences_water}--\ref{sec: alg} are put into action in Sections \ref{sec: water}--\ref{sec: result} by formulating and solving a water allocation problem and informing infrastructure decisions in Tucson, AZ.
Specifically, 
Section \ref{sec: water} provides the mathematical formulation and discusses how future water demands and supplies are predicted. 
Then, 
Section \ref{sec: result} presents the numerical results and analysis, and investigates the decentralized infrastructure decisions.
The paper ends in Section \ref{sec: concl_water} with 
concluding remarks.\vspace*{-0.09in}

\section{Multistage Distributionally Robust Optimization with $\phi$-Divergences}%
    \label{sec: phi_divergences_water}
\subsection{$\phi$-Divergences in DRO}
    \label{ssec: int_phi_divergences}
          Because we will be using $\phi$-divergences to form MDRO, we begin by reviewing $\phi$-divergences and presenting select properties of static/two-stage DROs 
        formed by $\phi$-divergences 
        that will be used later in the paper. 
 $\phi$-divergences provide a measure of distance between two distributions. 
    Let us focus on discrete distributions with a finite number of realizations  because we will work with finite scenario trees.
    Let $(\cdot)^T$ denote the transpose of a vector, and let $\mathbf{q} = (q_1, \dots, q_n)^T$ and $\mathbf{p} = (p_1, \dots, p_n)^T$ be two $n$-dimensional \textit{probability vectors}, i.e., satisfying $q_\omega, p_\omega \ge 0,$ for all $\omega=1,2,\ldots, n$ and $\sum_{\omega=1}^n p_\omega = \sum_{\omega=1}^n q_\omega = 1$.  
    The $\phi$-divergence from $\mathbf{p}$ to $\mathbf{q}$ is defined by
        \begin{equation}
            I_\phi(\mathbf{p},\mathbf{q}) = \sum_{\omega=1}^n q_\omega \phi\left(\frac{p_\omega}{q_\omega}\right),     \label{eq:I}
        \end{equation}
    where $\phi(u)$---called the {\it $\phi$-divergence function}---is a convex function on $u \geq 0$ such that $\phi(u) \geq 0$, $\phi(1) = 0$, and with the interpretations $0 \phi(a/0) = a \lim_{t \rightarrow \infty} \frac{\phi(u)}{u}$ and $0 \phi(0/0) = 0$.
    
    The $\phi$-divergence given in (\ref{eq:I}) is the expectation of the $\phi$-divergence function with respect to the nominal distribution $\mathbf{q}$, evaluated at the ratios $\frac{p_\omega}{q_\omega}$.
The convex conjugate of $\phi$ is defined as $\phi^*(s)=\sup_{u\geq 0}\{su-\phi(u)\}, s\in\mathbb{R}$.
A bound on the domain of the conjugate $\phi^*$ can be obtained through \textit{$\bar{s}:=\lim\limits_{u \rightarrow \infty} \frac{\phi(u)}{u}$}. That is, for any $s>\bar{s}$, $\phi^*(s)=\infty$.
The conjugate $\phi^*$ and the bound on its domain $\bar{s}$ will play an important role in reformulating the MDRO in Section \ref{ssec:  MDRO_reform} and the decomposition algorithm of Section \ref{sec: alg}. 
    
    Table \ref{tb:phi_examples} lists the $\phi$-divergences used for the water allocation problem, along with their conjugates.  %
    The modified $\chi^2$ distance is related to the famous $\chi^2$ goodness-of-fit test.
    Kullback-Leibler (KL) divergence is commonly used in probability and information theory.
    It can be interpreted as the expected log-scale loss $\sum p_\omega \left(\log\left({p_\omega}\right)-\log\left({q_\omega}\right)\right)$ with respect to $\mathbf{p}$.
    Hellinger distance is the squared Euclidean distance between $(\sqrt{p_1}, \dots, \sqrt{p_n})^T$ and $(\sqrt{q_1}, \dots, \sqrt{q_n})^T$. 
    Burg entropy changes the order of $\mathbf{p}$ and $\mathbf{q}$ in KL divergence; so it is the expected log-scale loss with respect to
 $\mathbf{q}$.

    \begin{table}[h]
    \small
        \centering
        \begin{tabular}{lcccc}
            \hline 
            Divergence                   & $\phi(u), u \geq 0$ & $\bar{s}$ &  $I_\phi(p,q)$ & $\phi^*(s)$ \\
            \hline
            Modified $\chi^2$ Distance   & $(u-1)^2$                     & $\infty$  & $\sum {(p_\omega - q_\omega)^2}/{q_\omega}$ & $\left\{\begin{array}{ll}
-1, & s < -2 \\
s+\frac{s^2}{4},& s\geq -2\end{array}\right.$
 \\
            KL Divergence  & $u\log u - u + 1$             & $\infty$  & $\sum p_\omega \log\left({p_\omega}/{q_\omega}\right)$ & $e^s-1$ \\
            Hellinger distance           & $\left(\sqrt{u}-1\right)^2$   & 1         & $\sum \left( \sqrt{q_\omega}-\sqrt{p_\omega}\right)^2$ & $\frac{s}{1-s}, s<1$\\
            Burg Entropy                 & $-\log u + u - 1$             & 1         & $\sum q_\omega \log\left({q_\omega}/{p_\omega}\right)$ & $-\log (1-s), s<1$\\ \hline
        \end{tabular}
                \caption{
            $\phi$-divergences used in this study.
        }
        \label{tb:phi_examples}
    \end{table}

Let us now briefly review DRO in the static or two-stage optimization context to reveal further properties of $\phi$-divergences used in this study.
Let $\xi$ be a random vector that takes values $\xi^1,\ldots,\xi^n$ with nominal probabilities $q_1,\ldots,q_n$.
DRO minimizes the worst-case expectation from a set of distributions that are similar---defined in a precise way below---to the nominal distribution $\mathbf{q}$. 
The resulting DRO formulation is
\begin{equation}
\min_{x \in \cX} \max_{\pbf \in \mathcal{P}} \ee{\pbf}{f(x,\xi)},\label{DRO}
\end{equation}
where the \textit{ambiguity set of distributions} is given by
$\mathcal{P} = \{\mathbf{p}: I_\phi(\mathbf{p},\mathbf{q}) \leq \rho,\ \sum_{\omega=1}^{n} p_\omega = 1,\ p_\omega \geq 0,\ \forall \omega \}$. Suppose for every $\xi$, $f(x,\xi)$ is a real-valued convex function on an open set containing $\cX$ and $\cX$ is a nonempty compact set. This ensures \eqref{DRO} has a finite optimal solution.
The first constraint in $\mathcal{P}$ only selects distributions sufficiently close to $\mathbf{q}$ with respect to a given $\phi$-divergence.
The remaining constraints in $\mathcal{P}$ ensure $\mathbf{p}$ is a probability vector.

The value of $\rho$ used in the first constraint in $\mathcal{P}$ determines the size of the ambiguity set. We refer to this parameter as \textit{the level of robustness.}
When $\phi$ is twice continuously differentiable around $1$ with $\phi^{\prime \prime}(1)>0$ (like those in Table \ref{tb:phi_examples}), 
$\rho$ can be defined as 
$\frac{\phi''(1)}{2N} \chi^2_{n-1,1-\alpha},$
where $N$ denotes the total number of observations 
and  $\chi^2_{n-1,1-\alpha}$ represents the $1-\alpha$ quantile of a chi-squared distribution with $n-1$ degrees of freedom. 
This value of $\rho$ produces an approximate $1-\alpha$ confidence region on the true distribution under mild conditions \citep{pardo2005statistical, bental2011robust}.

Let $(x^*,\pbf^*)$ be an optimal solution of \eqref{DRO}. We refer to $\pbf^*=(p_1^*,\ldots,p_n^*)^T$ as a \textit{worst-case probability vector} and $\ee{\pbf^*}{f(x^*,\xi)}$ as the \textit{worst-case expectation}. 
$\phi$-divergences differ in the way the worst-case probability vector can be formed. Suppose the nominal probability of scenario $\omega$ is positive, $q_\omega >0$. Some $\phi$-divergences are capable of \textit{suppressing} this scenario. That is, they may allow its worst-case probability to be zero, $p_\omega^*=0$. In essence, such a scenario is excluded from the final worst-case expectation. However, not all $\phi$-divergences are capable of suppression, and those that do, can suppress in different ways \citep{bayraksan2015data}. 

Among the $\phi$-divergences used in this study, modified $\chi^2$, KL divergence, and Hellinger distance are capable of suppressing scenarios.
Problem \eqref{DRO} formulated with the modified $\chi^2$ distance may choose to suppress any scenario individually, thus generating a wide variety of possible model output.
In contrast, when the KL divergence or the Hellinger distance is used, the only possible results are (a) no scenarios will be suppressed (i.e., $p^*_\omega>0$ for every $\omega$), or (b) all but the most costly scenarios will be suppressed.
Unlike the three $\phi$-divergences discussed above, the Burg entropy is not capable of suppressing scenarios. Thus, the solution will always have  $p^*_\omega>0$.
We will examine the implications of these behaviors in the multistage setting within the context of our application.
Next, we discuss how to extend the DRO problem \eqref{DRO} to the multistage setting, focusing on multistage linear optimization.\vspace*{-0.1in}

\subsection{MDRO Modeling and Formulation}
    \label{ssec: form_DRO}
In this paper we consider a discrete-time stochastic process in $T$ time stages, $\xi = (\xi_1, \xi_2, \ldots, \xi_T)$, where $\xi_t$ denotes the random vector composed of stochastic parameters ($A_t, B_t, b_t, c_t$) %
 of stage $t$  and $\xi_1$ is a degenerate random vector (i.e., a constant).  
We use $\xi_{[t]}=(\xi_1, \xi_2, \ldots \xi_t)$ to denote the history of the process through time $t$. 
We assume $\xi_t$ has a finite number of realizations for all stages $t=2,\ldots,T$; so the stochastic process can be represented as a finite scenario tree. 
We also assume the distribution governing the evolution of $\xi$ does not depend on the decisions. %
Note that we do not make any assumptions on the dependence structure of $\xi$. So, $\xi$ may be interstage independent, dependent according to a Markov structure, or can have more complicated dependencies.

Traditional multistage stochastic programs optimize a sequence of decisions at each stage $t$ that minimize the conditional expectation of an objective function, given the decision and history of the process up to that stage. 
The nested formulation of traditional multistage stochastic linear programs is given by \vspace*{-0.07in}
\begin{equation}
        \min_{x_{1} \in \cX_{1}} c_{1}x_{1} + \ee{\q_{2|\xi_{[1]}}}  {   \min_{x_{2} \in \cX_{2}(x_1, \xi_2)}  c_{2}x_{2} +\ee{\q_{3|\xi_{[2]}}} {  \ldots + \ee{\q_{T|\xi_{[T-1]}}} {   \min_{x_{T} \in \cX_{T}(x_{T-1}, \xi_T)} c_{T}x_{T} } \ldots } }, \label{msp} \vspace*{-0.07in}
\end{equation}
where $x_{t}:=x_{t}(\xi_{[t]})$ denotes the decisions (e.g., water allocations to different users, storage decisions at reservoirs, etc.) at stage $t$.
 The sequence of decisions $x_1, x_2, \ldots, x_T$ is collectively called a \textit{policy}.
Decisions $x_t$ only depend on the history of stochastic process up to stage $t$, i.e.\ $\xi_{[t]}$, and not the future.
This ensures  the decisions are \textit{nonanticipatory} and \textit{implementable}.
The multifunctions $\cX_{t}:=\cX_{t}(x_{t-1}, \xi_{t})=\{x_t:A_tx_t=B_tx_{t-1}+b_t, x_t\geq 0 \}$, for $t=2,\ldots,T$ and $\cX_1=\{x_1:A_1x_1=b_1, x_1\geq 0 \}$ 
represent the feasibility sets. 
We use the shorthand notation $x_{t}$ and $\cX_{t}$ to ease the presentation and switch to the full notation when we want to emphasize the dependencies.

The feasibility sets $\cX_{t}$ at stages $t=2,\ldots,T$ change according to the decisions of the previous stage 
$x_{t-1}(\xi_{[t-1]})$ and the stochastic parameters $\xi_t=(A_t, B_t, b_t, c_t)$ of that stage.
In the above formulation, $\q_{t|\xi_{[t-1]}}$ denotes the stage-$t$ conditional probability distribution, conditioned on the history of the process up to that point $\xi_{[t-1]}$, and $\ee{\q_{t}|\xi_{[t-1]}}{\cdot}$ denotes the conditional expectation taken with respect to  $\q_{t|\xi_{[t-1]}}$ for $t=2,\ldots, T$. 

The above model assumes the underlying probability distribution is known.
However, such an assumption is quite unrealistic, and typically a decision maker  only has partial information.
Our application particularly suffers from this issue, especially as it looks further into the future.
To address this limitation, a distributionally robust approach can be used. 

In this paper we build a distributionally robust problem in the multistage setting by constructing \textit{conditional} ambiguity set of distributions on a given scenario tree.
Note that the scenario tree can contain realizations of $\xi_t$ with zero conditional nominal probabilities.
However, we do not consider such realizations in our application. 
Also, at present, we suppress the notation for scenario trees to avoid cluttered exposition.
At each stage-$t$ ($t<T$) node of the tree, instead of using only one conditional distribution $\q_{t+1|\xi_{[t]}}$, an ambiguity set of conditional distributions are considered.
With this construction, the distributionally robust counterpart of \eqref{msp} is formulated as\vspace*{-0.09in} 
    \begin{equation}
    \begin{split}
    \min_{x_{1} \in \cX_{1}} c_{1}x_{1} + \max\limits_{\pbf_{2|\xi_{[1]}} \in \mathcal{P}_{2|\xi_{[1]}}}\mathbb{E}_{\pbf_{2|\xi_{[1]}}}   \left[   \min_{x_{2} \in \cX_{2}(x_1, \xi_2)}  c_{2}x_{2} + \max\limits_{\pbf_{3|\xi_{[2]}} \in \mathcal{P}_{3|\xi_{[2]}}}\mathbb{E}_{\pbf_{3|\xi_{[2]}}} \Bigg[ \ldots +  \Bigg. \right. \\
    \Bigg. \left. \max\limits_{\pbf_{T|\xi_{[T-1]}} \in \mathcal{P}_{T|\xi_{[T-1]}}}\mathbb{E}_{\pbf_{T|\xi_{[T-1]}}} \left[   \min_{x_{T} \in \cX_{T}(x_{T-1}, \xi_T)} c_{T}x_{T} \right]  \ldots  \Bigg]  \right] , \label{eq:nested_form}
    \end{split}
\end{equation}
where  $\mathcal{P}_{t+1|\xi_{[t]}}$ denotes the conditional ambiguity set, conditioned on the history of the stochastic process up to that stage $\xi_{[t]}$ for $t=1,\ldots, T-1$. 

There are various ways to construct the conditional ambiguity sets $\mathcal{P}_{t+1|\xi_{[t]}}$. 
As mentioned before, we focus on $\phi$-divergences. Similar to Section \ref{ssec: int_phi_divergences}, for $t=1,\ldots, T-1$,
we define the conditional ambiguity sets as\vspace*{-0.1in}
    \begin{align}
    \mathcal{P}_{t+1|\xi_{[t]}}=  \Big\{\pbf_{t+1|\xi_{[t]}}:  & \ \ \ I_{\phi}(\pbf_{t+1|\xi_{[t]}},\q_{t+1|\xi_{[t]}})\leq \rho_{t},  & \ \ \left(\lambda_{t}\right)  \label{eq:ambuguity_set}\\  
    & \ \ \ \ \mathbf{1}^T \pbf_{t+1|\xi_{[t]}} =1,&   \left(\mu_{t}\right) \nonumber \\
    & \  \ \ \ \ \ \ \ \pbf_{t+1|\xi_{[t]}} \geq 0 \Big\}, & \nonumber 
    \end{align}
    where $\mathbf{1}$ is a vector of the same size as $\pbf_{t+1|\xi_{[t]}}$ with all elements equal to 1.
    Although generally different $\phi$-divergences and different ambiguity sets can be used at various stages and histories of the process, in this paper,  we use the same $\phi$-divergence throughout the MDRO. 
    However, we change how close we remain to the nominal conditional distributions at different stages by changing the level of robustness $\rho_{t}.$
    We use the time index $t$ for $\rho_{t}$ and the dual variables ($\lambda_{t}$, $\mu_{t}$) of the first two constraints in \eqref{eq:ambuguity_set} because they will be used to reformulate the problems at stage $t$.
    
    Because $\phi$-divergences result in convex ambiguity sets, MDRO is a convex optimization problem.
   However, it becomes very difficult to solve as the number of stages $T$ and the number of realizations of $\{\xi_t\}^{T}_{t=1}$ increases.
   We will shortly present a formulation of MDRO that dualizes the inner maximization problems. This formulation will be utilized in our decomposition method.

    To ease the subsequent discussion, let us define stage-$t$ cost-to-go (value) functions as\vspace*{-0.03in} 
    \begin{align}
    \mathcal{Q}_{t}(x_{t-1}, \xi_{[t]})=\min_{x_{t} \in \cX_{t}(x_{t-1},\xi_t)} \  c_{t}x_{t} + \max_{\pbf_{t+1|\xi_{[t]}} \in \mathcal{P}_{t+1|\xi_{[t]}}} \ee{\pbf_{t+1|\xi_{[t]}}}{\mathcal{Q}_{t+1}(x_t, \xi_{[t+1]})} \label{eq:primal_form}
    \end{align}
    for $t=2,\ldots,T-1$. At the last stage $T$, the maximization problem in \eqref{eq:primal_form} is absent.
    At the first stage, we solve the following program\vspace*{-0.09in}
    \begin{align}
    \min_{x_{1} \in \cX_{1}} \  c_{1}x_{1} + \max_{\pbf_{2|\xi_{[1]}} \in \mathcal{P}_{2|\xi_{[1]}}} \ee{\pbf_{2|\xi_{[1]}}}{\mathcal{Q}_{2}(x_1, \xi_{[2]})}. \label{eq:primal_form_first_stage}
    \end{align}
    We assume $\cX_1 \neq \emptyset$, all feasibility sets $\cX_{t}$, $t=1,\ldots,T$ are bounded, and the problem \eqref{eq:primal_form_first_stage}--\eqref{eq:primal_form} has \textit{relatively complete recourse}. Together with the boundedness assumption this means that 
    for all $t=2,\ldots,T$, the feasibility sets $\cX_t$ are nonempty and bounded for any given feasible $x_{t-1}$ and any realization of $\xi_{[t]}$. Furthermore, we assume $\mathcal{Q}_{t}(x_{t-1}, \xi_{[t]})$ are finite for any given feasible $x_{t-1}$ and all realizations of $\xi_{[t]}$, for $t=2,\ldots,T$.\vspace*{-0.07in}
  
  \subsection{Risk-Averse Interpretation}
  \label{ssec:  MDRO_risk}
  Let us now discuss risk aversion in MDRO.
Recall that the so-called {\it risk measures} assign a value to each random outcome, indicating a preference between different outcomes of a random variable. \cite{Artzner1999}, in their pioneering work, argue that good risk measures should satisfy desirable properties like convexity and monotonicity, among others. They referred to such risk measures as {\it coherent risk measures}; see, e.g., \cite{shapiro2009} for further details.  Conditional Value-at-Risk (CVaR), for example, is one of the most popular coherent risk measures used today. We will now discuss relation of MDRO to such risk measures. 
  
  Consider problem \eqref{eq:primal_form} at stage $t=T-1$. 
  Observe that each conditional ambiguity set $\mathcal{P}_{T|\xi_{[T-1]}}$ is a bounded closed convex subset of (conditional) probability measures, which are defined over a finite set of realizations of $\xi_T$ given $\xi_{[T-1]}$.
  This and the fact that $\mathcal{Q}_{T}(x_{T-1},\cdot)$ is random, real-valued, i.e., a random variable, means the maximization problem in \eqref{eq:primal_form} forms (a conditional analogue of) a coherent risk measure; see, e.g, Theorem 3.1 of \cite{SHAPIRO2012719}. See also an alternative, axiomatic study of conditional risk mappings by \cite{RuszczynskiConditional2006}. Let us denote this risk measure as $\mathcal{R}_{\phi_{T|\xi_{[T-1]}}}$, where we suppress the dependence on $\rho_{T-1}$. So, $\mathcal{R}_{\phi_{T|\xi_{[T-1]}}}(\mathcal{Q}_{T}(x_{T-1}, \cdot))=\max_{\pbf_{T|\xi_{[T-1]}} \in \mathcal{P}_{T|\xi_{[T-1]}}}\ee{\pbf_{T|\xi_{[T-1]}}}{\mathcal{Q}_{T}(x_{T-1}, \cdot)}$. 
   Because of translation invariance property of coherent risk measures, we can equivalently write \eqref{eq:primal_form} at stage $T-1$ as $\mathcal{Q}_{T-1}(x_{T-2}, \xi_{[T-1]})=\min_{x_{T-1} \in \cX_{T-1}}\mathcal{R}_{\phi_{T|\xi_{[T-1]}}} \left(c_{T-1}x_{T-1}+\mathcal{Q}_{T}(x_{T-1}, \xi_{[T]})\right)$.
  Note that $\mathcal{Q}_{T-1}(x_{T-2}, \xi_{[T-1]})$ is finite for all possible values of $\xi_{[T-1]}$ and feasible $x_{T-2}$. Furthermore, it takes different values according to $\xi_{[T-1]}$ even when $x_{T-2}$ is fixed, and hence $\mathcal{Q}_{T-1}(x_{T-2}, \cdot)$ is a random variable. Going backward in time to $t=T-2,\ldots, 1$ and using the same reasoning as above recursively, we obtain a (conditional) coherent risk measure at each stage: $\mathcal{R}_{\phi_{t+1|\xi_{[t]}}}(\cdot)=\max_{\pbf_{t+1|\xi_{[t]}} \in \mathcal{P}_{t+1|\xi_{[t]}}}\ee{\pbf_{t+1|\xi_{[t]}}}{\cdot}$. We reach the following result.
  \begin{proposition}\label{pr:risk}
Consider the MDRO given in \eqref{eq:nested_form}, constructed via conditional ambiguity sets. 
Let $\mathcal{R}_{\phi_{t+1|\xi_{[t]}}}(\cdot)=\max_{\pbf_{t+1|\xi_{[t]}} \in \mathcal{P}_{t+1|\xi_{[t]}}}\ee{\pbf_{t+1|\xi_{[t]}}}{\cdot}$, where $\mathcal{P}_{t+1|\xi_{[t]}}$ is obtained through \eqref{eq:ambuguity_set} for $t=1,\ldots, T-1$. Then,  
\eqref{eq:nested_form} is equivalent to a multistage stochastic program with nested risk measures
    \begin{equation*}
    \begin{split}
    \min_{x_{1} \in \cX_{1}} c_{1}x_{1} + \mathcal{R}_{\phi_{2|\xi_{[1]}}}   \left[   \min_{x_{2} \in \cX_{2}(x_1, \xi_2)}  c_{2}x_{2} + \mathcal{R}_{\phi_{3|\xi_{[2]}}} \Bigg[ \ldots +  \Bigg. \right. 
    \Bigg. \left. \mathcal{R}_{\phi_{T|\xi_{[T-1]}}} \left[   \min_{x_{T} \in \cX_{T}(x_{T-1}, \xi_T)} c_{T}x_{T} \right]  \ldots  \Bigg]  \right]. \label{eq:nested_risk_averse}
    \end{split}
\end{equation*}
\end{proposition}

Consider composition of risk measures at stages $t$ and $t+1$,  $\mathcal{R}_{\phi_{t|\xi_{[t-1]}}}\circ\mathcal{R}_{\phi_{t+1|\xi_{[t]}}}(\cdot)=\mathcal{R}_{\phi_{t|\xi_{[t-1]}}}\left(\mathcal{R}_{\phi_{t+1|\xi_{[t]}}}(\cdot)\right)$.
Then the nested risk formulation of Proposition \ref{pr:risk} can be written using the composite risk measure $\overline{\mathcal{R}}=\mathcal{R}_{\phi_{2|\xi_{[1]}}}\circ\cdots\circ\mathcal{R}_{\phi_{T|\xi_{[T-1]}}}$.
Composite risk measures inherit many properties of their underlying risk measures, but they can be very complicated to write explicitly. Suppose each ambiguity set $\mathcal{P}_{t+1|\xi_{[t]}}$ is formed the same way, that is, using the same $\phi$-divergence and the same value of $\rho_t$. Even in this case, $\overline{\mathcal{R}}$ can be very different than any of the individual $\mathcal{R}_{\phi_{t+1|\xi_{[t]}}}$. One notable exception is the expectation risk measure, which satisfies the well-known equality $\e{\e{X|Y}}=\e{X}$ for any two random variables $X$ and $Y$. For further details on composite risk measures, we refer the readers to \cite{SHAPIRO2012719} and \cite{RuszczynskiConditional2006,RuszczynskiOptimization2006}. 

For many $\phi$-divergences, the explicit form of $\mathcal{R}_{\phi_{t|\xi_{[t-1]}}}$ is unknown, and as mentioned above the composite risk measure $\overline{\mathcal{R}}$ is even more complicated. Therefore, instead,  we will work with a dynamic-programming formulation of \eqref{eq:nested_form}  that is obtained through recursive application of Lagrangian duality. We present this formulation next.

        \subsection{MDRO Reformulation}    
    \label{ssec:  MDRO_reform}
    Suppose $\rho_t>0$ for $t=T-1,\ldots,1$. Then, the nominal conditional distribution $\q_{t+1|\xi_{[t]}}$ satisfies the first constraint in \eqref{eq:ambuguity_set} with $I_{\phi}(\q_{t+1|\xi_{[t]}},\q_{t+1|\xi_{[t]}})=0< \rho_{t}$; and so the Slater condition holds. Consequently, we have strong duality for the inner maximization problems in \eqref{eq:primal_form}--\eqref{eq:primal_form_first_stage}.
    Going backward from stage $T-1$ and dualizing the maximization problems in \eqref{eq:primal_form} with the Lagrangian multipliers  $\lambda_t$, $\mu_t$ from \eqref{eq:ambuguity_set}, we obtain \vspace*{-0.1in}
    \begin{alignat}{3}
    \mathcal{Q}_t(x_{t-1}, \xi_{[t]}) = &\min_{x_{t},\lambda_{t},\mu_{t}}  && c_{t}x_{t}+\ee{\q_{t+1|\xi_{[t]}}}{\mu_{t} +\rho_t\lambda_{t} + \lambda_{t} \phi^{*}\left( \frac{\mathcal{Q}_{t+1}(x_{t}, \xi_{[t+1]}) - \mu_{t}}{\lambda_{t}} \right)} \label{eq:dual_form}\\
    &\hspace{0.35cm}\mathrm{s.t.} && x_{t} \in \cX_{t}(x_{t-1},\xi_t), \ \lambda_{t} \geq 0, \nonumber\\
    & && \mathcal{Q}_{t+1}(x_{t}, \xi_{[t+1]}) - \mu_{t} \leq \bar{s}\lambda_{t}, \ \ \forall \xi_{t+1}|\xi_{[t]}, \nonumber
    \end{alignat}
    for $t=T-1, \ldots,2$,
    where $0\phi^*\left(\frac{a}{0}\right)=0$ when $a\leq 0$ and $0\phi^*\left(\frac{a}{0}\right)=+\infty$ if $a>0$.
    At the last stage, as before, $\mathcal{Q}_{T}(x_{T-1}, \xi_{[T]})=\min_{x_{T} \in \cX_{T}(x_{T-1},\xi_T)} \  c_{T}x_{T}$.
    The last constraint in \eqref{eq:dual_form} results from an \textit{implicit feasibility} consideration of the conjugate $\phi^*$
    based on its domain, which was discussed in Section \ref{ssec: int_phi_divergences}. 
    If $\bar{s}=\infty$ (see Table~\ref{tb:phi_examples}), this constraint is redundant and should be removed.
    Otherwise ($\bar{s}<\infty$), solutions that violate the last constraint cause the objective of \eqref{eq:dual_form} to be $\infty$; so they should not be considered.
    We explicitly present these constraints in the formulation because, for $\phi$-divergences with $\bar{s}<\infty$, our algorithm generates affine cutting planes to remove solutions that violate these constraints.
    At the first stage, we solve the following equivalent problem\vspace*{-0.05in}
     \begin{alignat}{3}
     &\min_{x_{1},\lambda_{1},\mu_{1}} && c_{1}x_{1}+\ee{\q_{2|\xi_{[1]}}}{\mu_{1} +\rho_1\lambda_{1} + \lambda_{1} \phi^{*}\left( \frac{\mathcal{Q}_{2}(x_{1}, \xi_{[2]}) - \mu_{1}}{\lambda_{1}} \right)} \label{eq:dual_form_first}\\
    &\hspace{0.35cm}\mathrm{s.t.} && x_{1} \in \cX_{1}, \ \lambda_1 \geq 0, \nonumber\\
    &  &&\mathcal{Q}_{2}(x_{1}, \xi_{[2]}) - \mu_{1} \leq \bar{s}\lambda_{1}, \ \ \forall \xi_{2}|\xi_{[1]}. \nonumber
    \end{alignat}
    
    \vspace*{-0.1in}
    \noindent
    Together, \eqref{eq:dual_form_first} and \eqref{eq:dual_form} provide a dynamic-programming formulation of MDRO when it is expressed as a minimization (instead of a minimax) problem.

    Such reformulations and their properties in the static/two-stage case are discussed in \cite{bental2011robust}; see also \cite{bayraksan2015data}. The above generalizes it to MDRO.
        When $T=2$---i.e., in the static/two-stage case---it is well known that DRO after dualization can be reformulated as a Second Order Cone Program (SOCP) when the modified $\chi^2$ or Hellinger distances are used 
    and as a convex program that admits a self-concordant barrier when the KL divergence or the Burg entropy are used.
    Such reformulations are referred to as \textit{robust counterpart problems}.
    A generalization to the MDRO follows the same: MDRO given in \eqref{eq:dual_form_first}--\eqref{eq:dual_form} is a large-scale SOCP if the conditional ambiguity sets are formed with the modified $\chi^2$ or Hellinger distances, and it is a large-scale convex program that admits a self-concordant barrier if the conditional ambiguity sets are formed with the KL divergence and the Burg entropy.
    As an example, we present the SOCP formulation\footnote{
    This formulation uses the scenario tree notation defined in Section~\ref{sec: alg}.}
    of MDRO formed with Hellinger distance in the Online Supplement.

    State-of-the-art solvers often take a long time to solve \eqref{eq:dual_form_first}--\eqref{eq:dual_form} with a large number of scenarios because the number of variables and constraints of the robust counterpart problem grows exponentially with the number of stages $T$.
    We next discuss a nested Benders decomposition algorithm to solve large-scale MDROs with $\phi$-divergences.\vspace*{-0.1in}

    \section{Decomposition Algorithm}
    \label{sec: alg}

    To explain the algorithm in a compact way, for $t=1,\ldots,T-1$, let $\tilde{\mathbf{x}}_t$ denote the collection of variables $(x_{t},\lambda_{t},\mu_{t})$ and let
        $\mathfrak{Q}_{t}(\tilde{\mathbf{x}}_{t}, \xi_{[t+1]})=\mu_{t} +\rho_t\lambda_{t} + \lambda_{t} \phi^{*}\left( \frac{\mathcal{Q}_{t+1}(x_{t}, \xi_{[t+1]}) - \mu_{t}}{\lambda_{t}} \right)$ represent the terms inside the expectations in \eqref{eq:dual_form_first}--\eqref{eq:dual_form}.
    The single-cut version of the  nested Benders algorithm replaces the convex functions $\ee{\q_{t+1|\xi_{[t]}}}{\mathfrak{Q}_{t}(\tilde{\mathbf{x}}_{t}, \xi_{[t+1]})}$ with a number of affine cutting planes to form their lower approximations for all $t=1,\ldots,T-1$. 
    When $\bar{s}<\infty$, the algorithm also generates affine feasibility cuts only when a candidate solution $\tilde{\mathbf{x}}_t$
    violates the implicit feasibility constraints.
    By removing the convex objective functions and the convex implicit feasibility constraints---and adding instead affine cuts for each---the algorithm solves only linear
    problems.
    
    To clearly present the algorithm, let us first describe our notation related to a scenario tree.
     The set of nodes at stage $t$ is denoted by $\Omega_t$, and $\omega_t\in\Omega_t$ denotes an element of this set, i.e., a stage-$t$ scenario. 
    By assumption, $\Omega_1$ is a singleton and $\xi_t$ has a finite sample space. Then, $\xi^{\omega_t}_t(\xi^{\omega_t}_{[t]})$ represents a specific realization of the random vector $\xi_t(\xi_{[t]})$. 
    A stage-$t$ ($t>1$) scenario $\omega_t$ has a unique ancestor in stage $t-1$, denoted by $a(\omega_t)$, and a stage-$t$ ($t<T$) scenario $\omega_t$ has a set of descendants, denoted by $\Delta(\omega_t)$. 
    To ease notation, we simply use $q^{\omega_{t}|\omega_{t-1}}$ instead of $q^{\omega_{t}|\omega_{t-1}}_{t|\xi_{[t-1]}} = P\big(\xi_{t}=\xi_{t}^{\omega_{t}}|\xi_{[t-1]}=\xi^{\omega_{t-1}}_{[t-1]}\big)$ to represent the nominal conditional probabilities on the nodes of the scenario tree. 
    Given this notation, all decision variables depend on $\omega_t$, i.e., $x^{\omega_t}_t$, $\lambda^{\omega_t}_t$, $\mu^{\omega_t}_t$, and $x_{t-1}$ is updated to $x^{a(\omega_t)}_{t-1}$.
    Earlier, we suppressed these dependencies for ease of exposition. 
    
    At node $\omega_t$ of stage $t$ ($t<T$), we have the following subproblem, denoted sub($\omega_t$): \vspace{-0.3cm}
    \begin{subequations}\label{eq:objective}
    \begin{alignat}{3}
    &\min_{\tilde{\mathbf{x}}_t^{\omega_t}, \theta^{\omega_t}_{t}} && c_t^{\omega_t}x_t^{\omega_t} + \theta^{\omega_t}_{t} &\label{eq:objective_func}\\
    & \hspace{0.4cm}\text{s.t.}  && A^{\omega_t}_{t} x^{\omega_t}_{t} = B^{\omega_t}_{t} x_{t-1}^{a(\omega_t)}+b^{\omega_t}_{t} 
    ,  & &\ \ \ \ \ (\pi^{\omega_t}_{t})\label{eq: setX}\\
    &  &&\theta^{\omega_t}_{t} \geq \sum_{\omega_{t+1}\in \Delta (\omega_t)} q^{\omega_{t+1}|\omega_t}\left(G^{\omega_{t+1}}_j  \tilde{\mathbf{x}}_t^{\omega_t} + g^{\omega_{t+1}}_j\right) , & \ \  \ \ \ j \in J^{\omega_t}_{t}, \label{eq:objective_cut}\\
    &  &&0\geq H^{\omega_{t+1}}_k  \tilde{\mathbf{x}}_t^{\omega_t} +h^{\omega_{t+1}}_k,  & \ \  \ \ \ k \in K^{\omega_t}_{t}, \label{eq:feasibility_cut}\\
    &  &&  x^{\omega_t}_t, \lambda^{\omega_t}_{t} \geq 0. &\nonumber
    \setlength\belowdisplayskip{0pt}
    \end{alignat} 
    \end{subequations}
    
    \vspace*{-0.1in}
    \noindent
    Constraints \eqref{eq:objective_cut} and \eqref{eq:feasibility_cut} represent the optimality and feasibility cuts, respectively. 
    Variable $\theta^{\omega_t}_{t}$, together with the optimality cuts \eqref{eq:objective_cut}, provide a lower approximation of $\ee{\q_{t+1|\xi_{[t]}}}{\mathfrak{Q}_{t}(\tilde{\mathbf{x}}_{t}, \xi_{[t+1]})}$. Feasibility cuts \eqref{eq:feasibility_cut} form an outer approximation of the implicit feasibility constraints.
    The sets $J^{\omega_t}_{t}$ and $K^{\omega_t}_{t}$ store the indices of all the cuts generated up to the current point in the algorithm.
    The subproblems at stage $T$ (sub($\omega_T$)) for all $\omega_T \in \Omega_T$ do not contain any cuts and do not have the decision variables $\theta_T^{\omega_T}, \mu_T^{\omega_T}$, $\lambda_T^{\omega_T}$; they only contain structural constraints (\ref{eq: setX}) and the non-negativity constraints $x^{\omega_T}_T\geq 0$.
    If $\bar{s}=\infty$, 
    there are no
    feasibility constraints \eqref{eq:feasibility_cut} %
    at any sub($\omega_t$).

    Let us now discuss how to obtain the cut coefficients, starting with the optimality cuts \eqref{eq:objective_cut}.
    Let ${\pi}^{\omega_{t}}_{t}$ denote the dual vector associated with the structural constraints \eqref{eq: setX}.
    At node $\omega_t\in\Omega_t$ of stage $t$ ($t<T$), suppose ($\hat{\tilde{\mathbf{x}}}_t^{\omega_t}$, $\hat{\theta}_t^{\omega_t}$) is a current solution of \eqref{eq:objective}. 
    We use $\hat{\cdot}$ to represent an optimal solution to \eqref{eq:objective} like $\hat{x}^{\omega_t}_{t}$, $\hat{\pi}^{\omega_t}_{t}$.
    To simplify the discussion, suppose $\hat{\lambda}_t^{\omega_t}>0$ and $\phi^*$ is differentiable (like those in Table \ref{tb:phi_examples}).
    When all descendant subproblems of $\omega_t$---that is, all sub($\omega_{t+1}$), $\omega_{t+1}\in\Delta(\omega_{t})$---are solved at $\hat{\tilde{\mathbf{x}}}_t^{\omega_t}$, we obtain the quantities\vspace*{-0.17cm}
    \[
    \hat{s}^{\omega_{t+1}}_{t+1} := ( c_{t+1}^{\omega_{t+1}}\hat{x}_{t+1}^{\omega_{t+1}}+\hat{\theta}^{\omega_{t+1}}_{t+1} -\hat{\mu}^{\omega_t}_{t})/{\hat{\lambda}^{\omega_t}_{t}},
    \]
    
    \vspace*{-0.13in}
    \noindent
    where the term $\hat{\theta}^{\omega_{t+1}}_{t+1}$ is absent when $t=T-1$.
    Then, the \textit{cut gradient} in the single-cut version of the algorithm is given by\vspace*{-0.1cm} $\sum\limits_{\omega_{t+1}\in\Delta(\omega_t)}q^{\omega_{t+1}|\omega_t}G^{\omega_{t+1}}_j$ and the \textit{cut intercept} is given by $\sum\limits_{\omega_{t+1}\in\Delta(\omega_t)}q^{\omega_{t+1}|\omega_t}g^{\omega_{t+1}}_j$, where\vspace*{-0.13in} 
    \begin{align*}
    G^{\omega_{t+1}}_j =& \Big(\phi^{*\prime}(\hat{s}^{\omega_{t+1}}_{t+1})  \hat{\pi}^{\omega_{t+1}}_{t+1} B^{\omega_{t+1}}_{t+1} \ \ \ \ \  \ 
    \rho^{\omega_t}_t + \phi^{*}(\hat{s}^{\omega_{t+1}}_{t+1}) -  \phi^{*\prime}(\hat{s}^{\omega_{t+1}}_{t+1})\hat{s}^{\omega_{t+1}}_{t+1} \ \ \  \ \  \ 
    1 - \phi^{*}{'}(\hat{s}^{\omega_{t+1}}_{t+1})\Big), \\
    g^{\omega_{t+1}}_j=&\hat{\mu}^{\omega_t}_t+\hat{\lambda}^{\omega_t}_t \rho^{\omega_t}_t+\hat{\lambda}^{\omega_t}_t\phi^{*}(\hat{s}^{\omega_{t+1}}_{t+1})-G^{\omega_{t+1}}_j \hat{\tilde{\mathbf{x}}}_t^{\omega_t}.
    \setlength\belowdisplayskip{0pt}
    \end{align*}
    
    \vspace*{-0.07in}
    The three terms of $G^{\omega_{t+1}}_j$ above correspond to the subgradients with respect to $x^{\omega_t}_t$, $\lambda^{\omega_t}_t$, and $\mu^{\omega_t}_t$, respectively.
    These quantities are obtained through the chain rule.  
    As an example, the first term of $G^{\omega_{t+1}}_j$ is calculated through
    $\frac{\partial \mathfrak{Q}_t}{\partial  x^{\omega_t}_{t}} = \phi^{*\prime}\left(\hat{s}^{\omega_{t+1}}_{t+1}\right)\cdot \frac{\partial\mathcal{Q}_{t+1}}{\partial  x^{\omega_t}_{t}} =  \phi^{*\prime}\left(\hat{s}^{\omega_{t+1}}_{t+1}\right)\cdot \hat{\pi}^{\omega_{t+1}}_{t+1} B^{\omega_{t+1}}_{t+1}$.
    The intercept term is obtained simply by using the subgradient inequality.
    The multi-cut version replaces $\theta^{\omega_t}_{t}$ in (\ref{eq:objective_func}) with $ \sum_{\omega_{t+1} \in \Delta (\omega_t)} q^{\omega_{t+1}|\omega_t}\theta^{\omega_{t+1}}_{t}$ and uses  individual cuts for each scenario $\omega_{t+1} \in \Delta(\omega_t)$ in (\ref{eq:objective_cut}):\  $\theta^{\omega_{t+1}}_{t} \geq G^{\omega_{t+1}}_j  \tilde{\mathbf{x}}_t^{\omega_t} + g^{\omega_{t+1}}_j$, for all $j \in J^{\omega_t}_{t}$. 
    We will compare the performance of the two variants in our numerical experiments.
    
    The cut coefficients of the feasibility cuts are obtained similarly. 
    If at $\omega_{t+1}\in \Delta(\omega_t)$, $\mathcal{Q}_{t+1}$ $\big(\hat{x}^{\omega_t}_{t}, \xi^{\omega_{t+1}}_{[t+1]}\big) - \hat{\mu}^{\omega_t}_{t} - \bar{s}\lambda_{t}>0$, we need to prevent solutions violating this constraint. 
    So, we must ensure the constraint is satisfied with $\leq 0$.
    By using subgradients and the chain rule, we obtain feasibility cuts \eqref{eq:feasibility_cut} with cut coefficients \vspace*{-0.1in}
    \begin{align*}
    H^{\omega_{t+1}}_k=&\Big(\hat{\pi}^{\omega_{t+1}}_{t+1} B^{\omega_{t+1}}_{t+1} \ \ \ \ -\bar{s} \ \ \ \  -1\Big),\\
    h^{\omega_{t+1}}_k=& c^{\omega_{t+1}}_{t+1}\hat{x}^{\omega_{t+1}}_{t+1}+\hat{\theta}^{\omega_t}_{t}-\hat{\pi}^{\omega_{t+1}}_{t+1}B^{\omega_{t+1}}_{t+1}\hat{x}^{\omega_t}_t.
    \setlength\belowdisplayskip{0pt}
    \end{align*}
    
    \vspace*{-0.1in}
    Algorithm \ref{alg:decomp} summarizes the single-cut version of the decomposition method.
    The algorithm works through two main phases: a \textit{forward pass} and a \textit{backward pass}.
    The forward pass solves all subproblems and stores all solutions.
    After the solutions' feasibility are checked and corrected, the algorithm updates the current upper bound. 
    The backward pass generates cutting planes to update the lower approximations.
    When the root node is solved with the current lower approximation, a lower bound to MDRO is obtained.
    Finally, the algorithm stops when the upper and lower bounds are sufficiently close.
    
    Algorithm \ref{alg:decomp} is a generalization of the so-called subgradient-based decomposition of \cite{NOYAN2012541}, originally developed for two-stage mean-CVaR stochastic programs. 
    Its extension to multistage mean-CVaR programs and its slight variant in \cite{shapiro_11} and \cite{kozmik_morton_15} are referred to as Decompositions \texttt{D4} and \texttt{D3}, respectively, in \cite{Weini2016}. 
    Observe that mean-CVaR multistage programs are equivalent to our setting when the $\phi$-divergence $\phi_{\text{EC}}$, defined as $\phi_{\text{EC}}(u)=0$ if $1-\kappa\leq u \leq 1+\frac{1}{1-\alpha}\kappa$ and $\phi_{\text{EC}}(u)=\infty$ otherwise, is used to form the conditional ambiguity sets in MDRO.
    This results in the conditional analogue of coherent risk measure 
    $\mathcal{R}_{\phi_{\text{EC}}}(\cdot)=(1-\kappa)\mathbb{E}[\cdot]+\kappa \text{CVaR}_\alpha(\cdot)$ at each stage.  Mean-CVaR multistage stochastic programs also arise when $L_\infty$-norm is used to form the conditional ambiguity sets \citep{Jianqiu2017Study}.  
    
    The mean-CVaR setting is significantly simpler than ours. 
    First, in its reformulation \eqref{eq:dual_form_first}--\eqref{eq:dual_form} there are no $\lambda_t$ variables. 
    Only the dual variables $\mu_t$ are present, and they represent the Value-at-Risk (VaR) in the CVaR representation $\text{CVaR}_\alpha(\cdot) = \min\limits_{\mu_t}\{\mu_t+\frac{1}{1-\alpha}\mathbb{E}[(\cdot-\mu_t)_+]\}$, where $(\cdot)_+$ denotes $\max\{0,\cdot\}$.
    Importantly, there are no implicit feasibility constraints in the mean-CVaR setting. 
    Algorithm \ref{alg:decomp} presents a significantly more general version of the subgradient-based decomposition, encompassing a large class of multistage distributionally robust and risk-averse multistage stochastic linear programs with nested coherent risk measures.
    
   We also note the decomposition algorithm of  \cite{PhilpottOnSolving2013}, further studied by \cite{Guigues2016}, and its specialization to modified $\chi^2$ distance \citep{philpott2018}. This algorithm, in addition to the linear subproblems at each node, explicitly solves the inner maximization problems in \eqref{eq:primal_form_first_stage}--\eqref{eq:primal_form}. Therefore, it solves additional convex programs at each stage-$t$ ($t<T$) node per iteration but always generates feasible solutions. 
   Algorithm \ref{alg:decomp} also always generates feasible solutions
   when $\bar{s}=\infty$ but only solves linear subproblems. The performance of different decomposition algorithms may depend on the underlying problem and specific instances.  Therefore, it is important to devise different decomposition algorithms, which can be successfully used and further specialized to specific problem classes and instances. 
   To the best of our knowledge, this is the first generalization of the subgradient-based algorithm to a large class of MDRO beyond the simple mean-CVaR setting.
    
        \begin{algorithm}[!th]
  \caption{Decomposition algorithm to solve MDRO}
\label{alg:decomp}
{\small
  \begin{algorithmic}
    \STEPZERO{Set $z_L = -\infty, z_U = +\infty$; Select small $\vartt{TOL}>0$, $\epsilon>0$; \hspace{10cm} 
                 For $t=1,\ldots,T-1$ and $\omega_t \in \Omega_t$, set  $J^{\omega_t}_{t}\leftarrow \emptyset, K^{\omega_t}_{t}\leftarrow \emptyset$; 
                 Add initial cuts; \\
                \label{alg:initial}}
    \item[]\vspace{0.1cm}
    \STEPSECOND{}\label{STEP2:s}
    \State{Solve sub($\omega_1$) and  obtain $\hat{x}^{\omega_1}_1$, $\hat{\lambda}^{\omega_1}_1$, $\hat{\mu}^{\omega_1}_1$, $\hat{\theta}^{\omega_1}_1$}
    \State{Set $z_L \leftarrow$ $c_1\hat{x}^{\omega_1}_{1} + \hat{\theta}_{1}$}
    \vspace{0.2cm} 
    \STEPFIRST{}\label{STEP1:s}
    \For{$t=2,\ldots,T$ and $\omega_t \in \Omega_t$}
    \State{\textbf{if} $t<T$ \textbf{then} Solve sub($\omega_t$) and  obtain $\hat{x}^{\omega_t}_t$, $\hat{\lambda}^{\omega_t}_t$, $\hat{\mu}^{\omega_t}_t$, $\hat{\theta}^{\omega_t}_t$ \textbf{end if}} 
    \State{\textbf{if} $t=T$ \textbf{then} Solve sub($\omega_T$) and obtain $\hat{x}^{\omega_T}_T$ and dual $\hat{\pi}^{\omega_T}_T$; Set $z^{\omega_T}_{T}\leftarrow c^{\omega_T}_{T}\hat{x}^{\omega_T}_{T}$\textbf{end if}}
    \EndFor\label{STEP1:e}
    \vspace{0.2cm} 
    \STEPTHIRD{}\label{STEP3:s}
    
    \For{$t=T-1,\ldots,1$ and $\omega_t \in \Omega_t$}
        \State{Set $\bar{\mu}_t^{\omega_t} \leftarrow \hat{\mu}_t^{\omega_t}$;}
        \State{\textbf{if} $(\hat{x}_t^{\omega_t}$, $\hat{\lambda}_t^{\omega_t}$, $\bar{\mu}_t^{\omega_t}$)  is infeasible according to \eqref{eq: infeasible_true} \textbf{then}  Find feasible $\bar{\mu}_t^{\omega_t}$  by using \eqref{eq: infeasible_true_adjust} \textbf{end if}}

            \State{Set $z^{\omega_t}_{t}\leftarrow c^{\omega_t}_{t}\hat{x}^{\omega_t}_{t}+\bar{\mu}^{\omega_t}_{t} +\rho^{\omega_t}_t\hat{\lambda}^{\omega_t}_{t} + \hat{\lambda}^{\omega_t}_{t} \sum\limits_{\omega_{t+1}\in\Delta(\omega_t)}q^{\omega_{t+1}|\omega_t}\phi^{*}\left( \frac{z^{\omega_{t+1}}_{t+1} - \bar{\mu}^{\omega_t}_{t}}{\hat{\lambda}^{\omega_t}_{t}} \right)$}
    \EndFor
        \If{$z^{\omega_1}_{1}< z_U$} 
           \State{Set $z_U \leftarrow z^{\omega_1}_{1}$; Update $x^{\omega_T*}_T \leftarrow \hat{x}^{\omega_T}_T$, \ $\forall \omega_T \in \Omega_T$;}
            \State{Update $(x^{\omega_t*}_t, \mu^{\omega_t*}_t, \lambda^{\omega_t*}_t)\leftarrow (\hat{x}^{\omega_t}_t, \bar{\mu}^{\omega_t}_t, \hat{\lambda}^{\omega_t}_t)$, $\forall \omega_t\in \Omega_t$, $\forall t<T$}
    \EndIf\label{STEP3:e}

    \vspace{0.2cm} 
    \STEPFOURTH{}\label{STEP4:s}
    \If{$z_U-z_L \leq \vartt{TOL}\cdot\min\left\lbrace\right |z_U|, |z_L|\rbrace$}
    \State{STOP}
    \State{Output: $x^{\omega_t*}_t, \ \forall \omega_t \in \Omega_t, \forall t$ with objective $z_U$ within $(100\cdot\vartt{TOL})\%$ of optimal }
    \EndIf\label{STEP4:e}
    \vspace{0.2cm} 
    \STEPFIFTH{}\label{STEP5:s}
    \For{$t=T-1,\ldots,1$ and $\omega_t \in \Omega_t$}
        \If{$(\hat{x}_t^{\omega_t}$, $\hat{\lambda}_t^{\omega_t}$, $\hat{\mu}_t^{\omega_t}$, $\hat{\theta}_t^{\omega_t})$  is infeasible according to \eqref{eq: infeasible_approx}}
            \State{Generate feasibility cut and add to problem (\ref{eq:objective}); Update $K^{\omega_t}_{t}$ \label{alg:gen_cut1}}
            \State{Adjust $\hat{\mu}^{\omega_t}_{t}$ to obtain a feasible $(\hat{x}_t^{\omega_t}$, $\hat{\lambda}_t^{\omega_t}$, $\hat{\mu}_t^{\omega_t}$, $\hat{\theta}_t^{\omega_t})$ by using \eqref{eq: infeasible_approx_adjust}\label{alg:adj_mu}}
        \EndIf
        \State{Generate objective cut and add to problem (\ref{eq:objective}); Update $J^{\omega_t}_{t}$; \label{alg:gen_cut2}}
        \State{\textbf{if} $t>1$ \textbf{then}  Solve sub($\omega_t$) to obtain $\hat{x}^{\omega_t}_t, \hat{\lambda}^{\omega_t}_t,\hat{\mu}^{\omega_t}_t, \hat{\pi}^{\omega_t}_t, \hat{\theta}^{\omega_t}_{t}$ \textbf{end if}\label{STEP5:e}}
    \EndFor
    \State{Go to STEP 1.}
  \end{algorithmic}
}
\end{algorithm}

    Let us now discuss select features of Algorithm \ref{alg:decomp}, focusing especially on feasibility because this feature is not present in simpler forms of the subgradient-based decomposition. 
    Algorithm \ref{alg:decomp} checks feasibility at two points: upper bound calculation (Step 3) and backward pass (Step 5). 
    We discuss these in this order next.
    To ease exposition, we assume all $\lambda^{\omega_t}_{t}>0$ in the subsequent discussion.
    
    The upper bound, calculated in Step 3, is  the objective function value of problem \eqref{eq:dual_form_first} with a feasible policy $\hat{{x}}_1, \hat{{x}}_2, \ldots, \hat{{x}}_T$.
    Because we don't keep all the implicit feasibility constraints, 
    there is a possibility of having an infeasible solution. 
    To test feasibility, for $t=1,\ldots, T-1$ and $\omega_t \in \Omega_t$, let us 
    keep track of 
    $\bar{\mu}^{\omega_t}_{t}=\hat{\mu}^{\omega_t}_{t}$,  
    which may change if infeasibility is detected,
    and  set
    $z^{\omega_t}_{t}=c^{\omega_t}_{t}\hat{x}^{\omega_t}_{t}+\bar{\mu}^{\omega_t}_{t} +\rho^{\omega_t}_t\hat{\lambda}^{\omega_t}_{t} + \hat{\lambda}^{\omega_t}_{t} \sum_{\omega_{t+1}\in\Delta(\omega_t)}q^{\omega_{t+1}|\omega_t}\phi^{*}\left( \frac{z^{\omega_{t+1}}_{t+1} - \bar{\mu}^{\omega_t}_{t}}{\hat{\lambda}^{\omega_t}_{t}} \right)$.
    At the last stage, $z^{\omega_T}_T$ is simply given by $z^{\omega_T}_{T}=c^{\omega_T}_{T}\hat{x}^{\omega_T}_{T}$. 
    A current solution is infeasible  with respect to \eqref{eq:dual_form_first}--\eqref{eq:dual_form} if at any $t=1, \ldots, T-1$\vspace*{-0.1in}
    \begin{equation}\label{eq: infeasible_true}
        \bar{s} < \infty \ \ \ \text{and} \ \ \ \sup\limits_{\omega_{t+1} \in \Delta(\omega_t)} 
        \tfrac{z^{\omega_{t+1}}_{t+1}-\bar{\mu}^{\omega_t}_{t}}
        {\hat{\lambda}^{\omega_t}_{t}} 
        > \bar{s}. %
    \end{equation}
    
    \vspace*{-0.07in}
    \noindent
    If the solution is infeasible, we adjust $\bar{\mu}^{\omega_t}_{t}$ to be feasible with the equation \vspace{-0.07in}
    \begin{eqnarray}\label{eq: infeasible_true_adjust}
    \bar{\mu}^{\omega_t}_{t}\leftarrow \sup\limits_{\omega_{t+1}\in \Delta(\omega_t)} z^{\omega_{t+1}}_{t+1}-\bar{s}\hat{\lambda}^{\omega_t}_t \left(1-\epsilon\right)  \vspace{-0.05in}
    \end{eqnarray}
    
    \vspace*{-0.07in}
    \noindent
    for some small $\epsilon>0$.
    A potential upper bound is calculated with this adjusted feasible solution, and the upper bound is updated if this value is smaller than the current upper bound. 
    
    The backward pass (Step 5) generates feasibility and optimality cuts going backward from stage $t=T-1$ down to $1$.
    If the current solution is infeasible to a subproblem \eqref{eq:objective}, i.e.,\vspace*{-0.17in}
    \begin{eqnarray}\label{eq: infeasible_approx}
    \bar{s} < \infty \ \ \  \text{and} \ \ \ \sup\limits_{\omega_{t+1} \in \Delta(\omega_t)} \Hat{s}^{\omega_{t+1}}_{t+1} >\bar{s},
     \setlength\belowdisplayskip{0pt}
    \end{eqnarray}
    
    \vspace*{-0.07in}
    \noindent 
    we generate a feasibility cut.
    Then, we adjust $\hat{\mu}^{\omega_t}_{t}$ to obtain a feasible solution
    with the below simple equation with a small $\epsilon>0$ and by setting  $\hat{\theta}^{\omega_{T}}_{T}\equiv 0$:\vspace*{-0.1in}
    \begin{equation}\label{eq: infeasible_approx_adjust}
    \hat{\mu}^{\omega_t}_{t}\leftarrow \sup\limits_{\omega_{t+1}\in \Delta(\omega_t)}c_{t+1}^{\omega_{t+1}}\hat{x}_{t+1}^{\omega_{t+1}} + \hat{\theta}^{\omega_{t+1}}_{t+1} -\bar{s} \hat{\lambda}^{\omega_t}_t \left(1-\epsilon\right).
    \end{equation}
    
    \vspace*{-0.1in}
    \noindent
    The adjusted feasible $\hat{\mu}^{\omega_t}_{t}$ is used to generate an optimality cut.
    
    At the beginning of backward pass 
    (i.e., when $t=T-1$), feasibility conditions \eqref{eq: infeasible_approx} and \eqref{eq: infeasible_true} and the adjustments \eqref{eq: infeasible_approx_adjust} and \eqref{eq: infeasible_true_adjust} are equivalent.  
    Once feasibility and optimality cuts are added to the subproblems in the backward pass, the subproblems are solved again, and an improved approximation to stage-$t$ problem is obtained. Starting from $T-2$ (for any $t<T-1$), a feasibility cut is added if a candidate solution is infeasible with respect to subproblem \eqref{eq:objective}, which is an \textit{approximation} of the original stage-$t$ problem; see condition \eqref{eq: infeasible_approx}. This cut is valid for the original problem because an infeasible solution of  \eqref{eq:objective} is also infeasible to the original problem. 
    In contrast, we ensure feasibility with respect to the \textit{original problem} \eqref{eq:dual_form_first}--\eqref{eq:dual_form} to calculate an upper bound in Step 3; see condition \eqref{eq: infeasible_true}.
    This is because we need a feasible solution to the original problem to calculate an upper bound.

    In implementation, we treat ${\lambda}^{\omega_t}_t>0$ and test optimality and feasibility at $\lambda^{\omega_t}=0$ separately.
    For example, in our water allocation problem, we set ${\lambda}^{\omega_t}_t\geq $ 1E-5 for $t<T$ and $\omega_t\in\Omega_t$, and we check $\lambda^{\omega_t}_{t}=0$ if ${\lambda}^{\omega_t}_t$ hits the lower bound of 1E-5.
    Furthermore, at the first iteration, for $t=1,\ldots, T-1$ and $\omega_t\in\Omega_t$, each sub($\omega_t$) needs an initial cut---e.g., $\theta^{\omega_t}_{t}\geq -M$ with large enough $M\geq 0$---to avoid unboundedness.
    In our numerical experiments, we set $\theta^{\omega_t}_{t}\geq 0$ as initial cuts for all $t<T$ and $\omega_t\in\Omega_t$ because  all subproblems have a lowest cost of zero in our application.
    
    Finally, we remark that because Algorithm \ref{alg:decomp} is based on the reformulated MDRO in \eqref{eq:dual_form_first}--\eqref{eq:dual_form}, it does not directly calculate the worst-case probabilities present in the original formulation \eqref{eq:nested_form} of MDRO. 
    Given $(\hat{x}^{\omega_t}_{t},\hat{\mu}^{\omega_t}_{t}, \hat{\lambda}^{\omega_t}_{t})$, the conditional probabilities $p^{\omega_{t+1}|\omega_{t}}$ can be easily calculated with the below three equations:\vspace*{-0.19in}
    \begin{eqnarray}\label{eq:p_worst}
        &&p^{\omega_{t+1}|\omega_{t}} = q^{\omega_{t+1}|\omega_t}\cdot\phi^{*'}\left(\frac{\mathcal{Q}_{t+1}\big(\hat{x}^{\omega_t}_{t}, \xi^{\omega_{t+1}}_{[t+1]}\big) - \hat{\mu}^{\omega_t}_{t}}{\hat{\lambda}^{\omega_t}_{t}}\right), \\
        &&\sum_{\omega=1}^n q^{\omega_{t+1}|\omega_t}\phi\left(\frac{p^{\omega_{t+1}|\omega_{t}}}{q^{\omega_{t+1}|\omega_t}}\right) = \rho^{\omega_t}_{t},\ \ \ \sum_{\omega_{t+1}\in \Delta (\omega_{t})}p^{\omega_{t+1}|\omega_{t}} = 1.\nonumber
    \end{eqnarray} 
    
    \vspace*{-0.1in}
    \noindent
    Furthermore, given $\vartt{TOL}'>0$, a secondary stopping rule  $(1-\sum_{\omega_{t+1}\in\Delta(\omega_t)}p^{\omega_{t+1}|\omega_{t}})\leq \vartt{TOL}'$, for all $t<T$, $\omega_t\in\Omega_t$ (or only at $t=1$) can be used in Step 4. 
    Because we want to examine the worst-case probabilities assigned to different climate models and demand scenarios, we also use this additional stopping rule in our implementation. 
    \vspace*{-0.1in}

\section{Application to Water Allocation Problem}
	\label{sec: water}\vspace*{-0.05in}
	\subsection{Problem Description}
	\label{ssec: probelm_des}
	\vspace*{-0.07in}
	
\begin{figure}[ht]
		\centering
		\includegraphics[width=30pc]{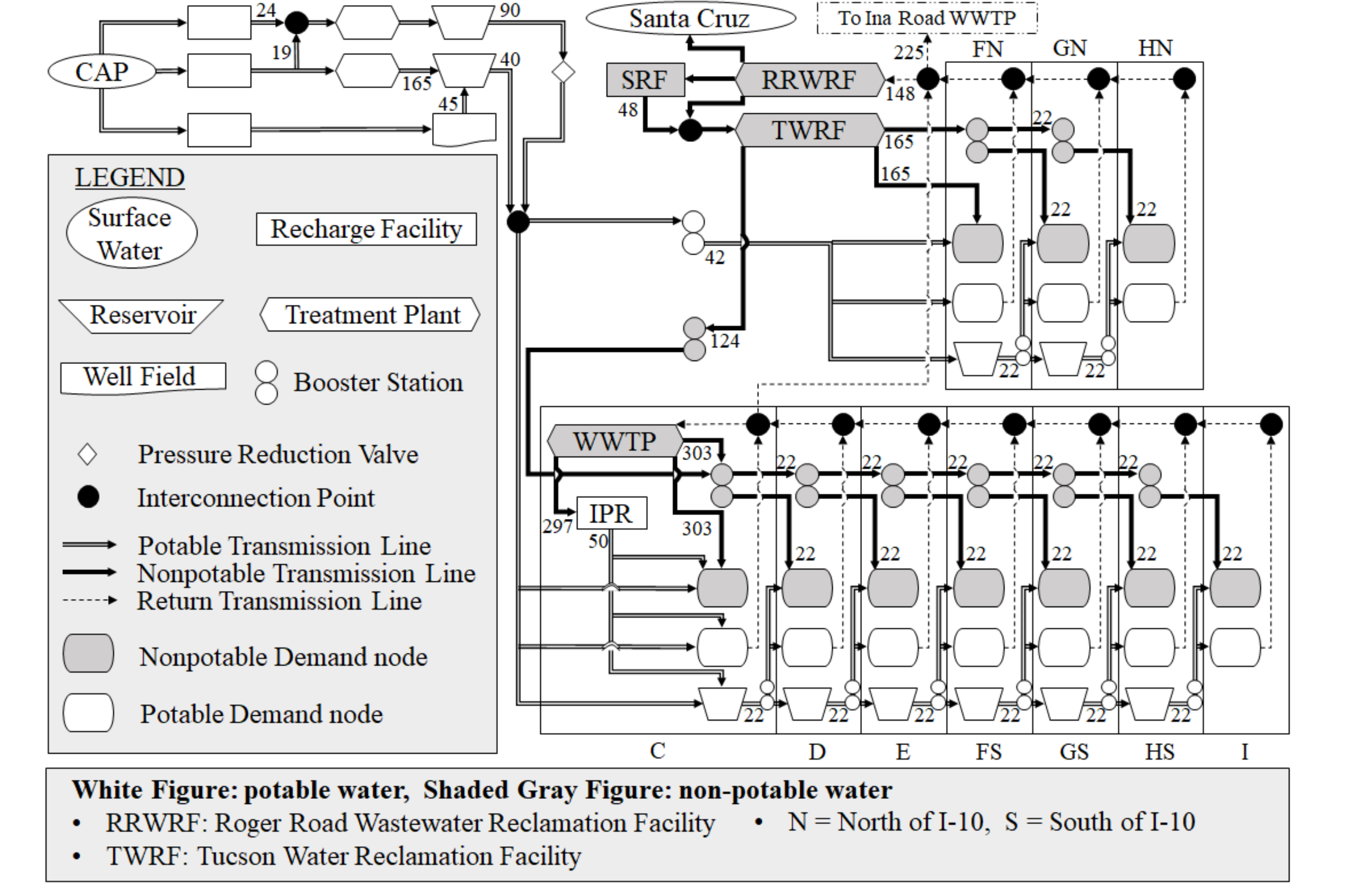}
		\caption{A schematic of the water system in the study area.}
		\label{fig:tucson_treatment}
	\end{figure}

    The southeastern portion of Tucson---called the \textit{study area}---is being increasingly developed.
	A schematic view of the area's water system is shown in Figure \ref{fig:tucson_treatment}.
	Majority of Tucson's water comes from the Colorado River, brought in by the Central Arizona Project (CAP) canal.
	This water is then treated and sent to customers, or seeped into underground to be saved for future use.
	These are represented as ``CAP'' and a few other white nodes in top-left corner of Figure \ref{fig:tucson_treatment}.

	This area is split into different demand zones: C, D, E, FS, FN, $\ldots$, I. 
	Given the limited capacity of the existing treatment plants and increasing population, the governing agencies in Tucson are interested in building additional treatment facilities in the area.
	Figure \ref{fig:tucson_treatment} shows both existing infrastructures (e.g., CAP) and proposed new infrastructures---a satellite wastewater treatment plant (WWTP) and Indirect Potable Reuse (IPR) facility---in Zone C.  
		The schematic also shows both the potable water system (through white nodes and double-lined acrs) and the reclaimed water system (through gray nodes and solid black arcs).	
		The wastewater return pipes are shown in dashed lines.\vspace*{-0.1in}
	
	Each zone contains potable and nonpotable demand nodes and a reservoir and booster station for transporting each type of water.
	Potable water, being of higher quality, can satisfy either type of demand.
	Figure \ref{fig:tucson_treatment} provides  the cost in \$/acre-foot (af) on each arc, if it is not negligible. 
	A dummy node capable of supplying water in the event of a water shortage is also included in the model (but not shown in figure).
	The cost of this extra supply is set at \$800/af by a fixed contract.
	In other words, there exists a water market at a constant exogenous price, as commonly used in the literature \citep{Murali2015,Calatrava2005,Weinberg1993}.
	
	The model aims to allocate the Colorado River water to different users in the study area while being sustainable through the mid-century. 
	Furthermore, the model is used to evaluate several local infrastructure decisions by water authorities who are facing considerable uncertainties.
	We present its formulation next.\vspace*{-0.1in}

    \subsection{Formulation}
    \label{ssec: form_water}
    The water system is represented as a directed network graph $\mathcal{G} = (\mathcal{N}, \mathcal{A})$, where $\mathcal{N}$ is the set of nodes and $\mathcal{A}$ is the set of arcs.
	For each year, the network has $|\mathcal{N}|=60$ nodes, 
	categorized into seven sets:
    pumps and reservoirs (PR), 
    water treatment plants (TP), 
	potable users (PU), 
	nonpotable users (NU), 
	recharge facilities (RF),
	water supply from the Colorado River (CAP), and
	a dummy node (D).
	The network also has $|\mathcal{A}|=100$ arcs, representing the pipes carrying water and 
	connecting the network to the five reservoirs for water storage to be used in future periods.

	We study a total of $P = 33$ time periods with 5 stages, representing years 2018--2050.
	Each stage $t$ comprises of several periods $p=1, \ldots, P_t$, where $P_1=1$ and $P_i=8$ for $i=2,\ldots,5$.
    Costs for all time periods are brought into present value by applying a 4\% discount rate per year.
    
	For each scenario $\omega_t \in \Omega_t$ of stage $t$, water flows on arc $(i,j) \in \mathcal{A}$ during time period $p = 1, \dots, P_t$ are represented by decisions $x^{\omega_t}_{i,j,p,t}$.
	Each arc $(i,j) \in \mathcal{A}$ at time period $p$ has a unit cost $c_{i,j,p,t}$ and loss coefficient $0 \leq a_{i,j,p,t} \leq 1$ to account for evaporation and leakage from the pipes.
	Stored water available at node $i$ at the end of time period $p$ of stage $t$ is denoted $x^{s,\omega_t}_{i,p,t}$, with associated storing cost $c_{i,p,t}^{s}$.
	Water released into the environment from node $i$ in period $p$ is similarly represented by $x^{r,\omega_t}_{i,p,t}$.
	We assume release into the enviroment has no cost.
	Water shortage decisions are represented by $x^{short,\omega_t}_{i,j,p,t}$ with associated shortage cost $c_{i,j,p,t}^{short}$.
    So, the decision vector consists of  $x^{\omega_t}_t=\left(x^{\omega_t}_{i,j,p,t},x^{s,\omega_t}_{i,p,t},x^{r,\omega_t}_{i,p,t},x^{short,\omega_t}_{i,j,p,t}\right)$.
	Also, we define capacities of the following nodes: recharge facility storage capacity ($U^{RF}_{i,p,t}$),
	pumping capacity at treatment plant or recharge facility ($U^{RFTP}_{i,p,t}$), 
	and treatment storage capacity ($U^{TP}_{i,p,t}$). 
	
	For a realization $\xi^{\omega_t}_{p,t}$ at period $p$ of stage $t$, $d_{i}(\xi^{\omega_t}_{p,t})$ denotes the demand of user node $i$ and $CAP(\xi^{\omega_t}_{p,t})$ represents the Colorado River water allotment to the study area. Notation $l$ denotes the portion of potable municipal user demand returned to a wastewater treatment facility.
	To compute the first-year constraints, we assume initial storage levels $x^{s}_{i,P_0,0}$ are given at recharge facilities. 

    Stage-$t$ ($t<T$) minimax formulation at node $\omega_t \in \Omega_t$ is given by\vspace*{-0.1in}
    \begin{eqnarray*}
    \mathcal{Q}\left(x^{s,a(\omega_t)}_{i,P_{t-1},t-1},\xi^{\omega_t}_{[t]}\right) = \min_{x^{\omega_t}_t\in \cX^{\omega_t}_t}\hspace{-0.5cm}&&\sum_{p=1}^{P_t} \Bigg(\sum_{(i,j) \in A}  c_{i,j,p,t} x^{\omega_t}_{i,j,p,t} + \sum_{(D,j) \in A}  c^{short}_{D,j,p,t} x^{short,\omega_t}_{D,j,p,t}+\sum_{j \in N}  c_{j,p,t}^{s} x^{s,\omega_t}_{j,p,t}\Bigg)\\
    &&+ \max_{\pbf^{\omega_t}_{t+1} \in \mathcal{P}^{\omega_t}_{t+1|\xi_{[t]}}}\sum_{\forall \omega_{t+1} \in \Delta (\omega_t)} p^{\omega_{t+1}|\omega_t} \mathcal{Q}_{t+1}(x^{\omega_t}_t, \xi^{\omega_{t+1}}_{t+1}),
    \end{eqnarray*}
    
    \vspace*{-0.1in}
    \noindent 
    where $\cX^{\omega_t}_t$ comprises of the following constraints\vspace*{-0.1in} 
	\begin{alignat}{3}
    &\sum_{j : (j,i) \in A} a_{j,i,p,t} x^{\omega_t}_{j,i,p,t} = \sum_{j : (i,j) \in A} x^{\omega_t}_{i,j,p,t}, & & \hspace{.3cm}i \in PR \cup TP, \ 1 \leq t \leq P_t,  \label{eq:ha}\\
	&\sum_{j : (j,i) \in A} a_{j,i,p,t} x^{\omega_t}_{j,i,p,t} + x^{short,\omega_t}_{D,i,p,t} = d_{i}(\xi^{\omega_t}_{p,t}), & & \hspace{.3cm}i \in PU \cup NU, \ 1 \leq t \leq P_t,  \label{eq:hb} \\
	& \sum_{j : (j,i) \in A} a_{j,i,1,t} x^{\omega_t}_{j,i,1,t} + x^{s,a(\omega_t)}_{i,P_{t-1},t-1}= \sum_{j : (i,j) \in A} x^{\omega_t}_{i,j,1,t}  + x_{i,1,t}^{r,\omega_t}+x_{i,1,t}^{s,\omega_t},& & \hspace{.3cm}i \in RF,   \label{eq:hc}\\
	& \sum_{j : (j,i) \in A} a_{j,i,p,t} x^{\omega_t}_{j,i,p,t} + x^{s,\omega_t}_{i,p-1,t}= \sum_{j : (i,j) \in A} x^{\omega_t}_{i,j,p,t} + x_{i,p,t}^{r,\omega_t}+ x_{i,p,t}^{s,\omega_t},  & &\hspace{.3cm}i \in RF,  \  2 \leq p \leq P_t,  \label{eq:hd}\\
	&   \sum_{j : (i,j) \in A}  x^{\omega_t}_{i,j,1,t} \leq x^{s,\omega_t}_{i,P_{t-1},t-1},& & \hspace{.3cm}i \in RF,      \label{eq:he}\\
	&   \sum_{j : (i,j) \in A}  x^{\omega_t}_{i,j,p,t} \leq x^{s,\omega_t}_{i,p-1,t},& & \hspace{.3cm}i \in RF, \  2 \leq p \leq P_t,   \label{eq:hf}\\
	&   \sum_{j : (i,j) \in A}  x^{\omega_t}_{i,j,p,t} \leq CAP(\xi^{\omega_t}_{p,t}),& &   \hspace{.3cm} 1 \leq p \leq P_t,   \label{eq:hg}\\
	&   \sum_{j : (i,j) \in A}  x^{\omega_t}_{i,j,p,t} \leq U^{RFTP}_{i,p,t},& & \hspace{.3cm}i \in RF \cup TP, \   1 \leq p \leq P_t,   \label{eq:hh}\\
	&   \sum_{j : (i,j) \in A}  x^{\omega_t}_{j,i,p,t} \leq U^{TP}_{i,p},& & \hspace{.3cm}i \in TP,  \  1 \leq p \leq P_t,   \label{eq:hi}\\
	& x^{\omega_t}_{i,j,p,t} = l\cdot d_{i}(\xi^{\omega_t}_{p,t}),&  &\hspace{-.3cm}i \in PU, (i,j) \in \mathcal{A},\   1 \leq p \leq P_t,   \label{eq:hj}\\
	&  0 \leq x_{i,p,t}^{s,\omega_t} \leq U^{RF}_{i,p,t},& & \hspace{.3cm}i \in RF,  \  1 \leq p \leq P_t,   \label{eq:hk}\\
	& x^{\omega_t}_{i,j,p,t} \geq 0,& & \hspace{.3cm}j : (i,j) \in \mathcal{A} ,  \  1 \leq p \leq P_t. \label{eq:hm}  
	\end{alignat}
	The constraints can be summarized into three categories. First, flow balance constraints on nodes include 
	(i) water flow balance at pumps/water treatment plants/reservoirs/interconnection points (\ref{eq:ha});
	(ii) demand satisfaction at potable/nonpotable users (\ref{eq:hb}); 
	(iii) water storage balance at recharge facilities \{(\ref{eq:hc}), (\ref{eq:hd})\}, where infiltration needs a one-year lag \{(\ref{eq:he}), (\ref{eq:hf})\}. Second, capacity constraints on nodes include
	(i) bounds on the Colorado River water supply depending on scenario (\ref{eq:hg}); and 
	(ii) bounds on the in/out-flow of recharge facilities and treatment plants \{(\ref{eq:hh}), (\ref{eq:hi})\}. Finally, constraints regarding arcs entail
	(i)  a fixed portion of the potable used water is returned to a wastewater treatment plant (\ref{eq:hj})---the treated water can be used only for nonpotable demand for later years---and (ii) upper bound and non-negativity constraints on the water flows \{(\ref{eq:hk}), (\ref{eq:hm})\}. 

The above model is similar to the one in \cite{Weini2016}. However, the stochastic representation of  \cite{Weini2016}'s model is quite crude. The model of this paper, in contrast, has significantly more detailed and realistic scenarios. For instance, it considers climate uncertainty for the first time. It has a substantially more detailed portrayal of the Colorado River water availability based on hydrological studies. Per-capita demand models and population estimates are also considerably improved by using various statistical methodologies and studies conducted by the local governing agencies. Furthermore, this model is based on MDRO with general $\phi$-divergences, which is a large generalization of the simpler mean-CVaR model in \cite{Weini2016}. Finally, there are no infrastructure decisions considered in that paper, whereas here, we use the MDRO model to evaluate these important decisions on constructing decentralized infrastructures.\vspace*{-0.1in}  
	\subsection{Scenario Generation}\vspace*{-0.1in}  
	\label{sec: ScenarioGeneration}
	Our optimization model requires two primary uncertain data:  annual water demand by zone and annual water supply (e.g, right-hand sides of (\ref{eq:hb}) and (\ref{eq:hg})).
	To quantify these, we use a large number of data from various sources, some developed by experts  in their fields. We summarize our data and its sources in the Online Supplement and highlight some important sources and our methodology below. \vspace*{-0.17in}

	\subsubsection{Annual Population Estimates}\vspace*{-0.07in}
	\label{ssec: population}

		Because population affects both water demand and supply, we discuss it first. 
	Water demand in each zone is \textit{proportional} to the population of each zone. 
	The Colorado River water allocation (=water supply),
	on the other hand, depends on the {\it ratio} of the study area's population to the overall Tucson population.  
	Therefore, we need population estimates for both the study area and Tucson.
	We used various local studies for 2050 population predictions for both the study area and Tucson. Then, we utilized the last U.S.\ census  numbers  \citep{uscensus} and interpolated the intermediate years. 
	The total population in the study area is then broken down to demand zones for each year, based on the population propagation model developed by Tucson Water and the City of Tucson. 
	The beginning (2010 census) and ending (2050 estimates) total population numbers and their sources are summarized in Table O.1 in the Online Supplement. As a result of this analysis, we have two population estimates in our model: 
	(i) a low-population and 
	(ii) a high-population estimate.

	\subsubsection{Water Demand Prediction}
	\label{ssec: demand_prediction}
	We first investigated how climate variables like temperature and precipitation as well as water-use trends affect water demand---measured in Gallons Per Capita per Day (GPCD)---by using historical data. 
	The historical data reveals that the average GPCD began dropping near the beginning of the 21st century, from over 170 GPCD to 140.
	This analysis produced two regression models: one that assumes increasing efficiency in water use (called the lower-GPCD model) and the other not (called the higher-GPCD model). 
	The lower-GPCD model might be appropriate if technological advances and water conservation efforts lead to significantly lower water consumption in the future. The higher-GPCD model, on the other hand, assumes people cannot decrease water consumption indefinitely.
	These models are then used to predict future demands by incorporating climate predictions with a given greenhouse concentration pathway. 
	We summarize this analysis below.
	\medskip 
	
	\noindent 
	\textbf{Building Statistical Models.} Because residuals of the ordinary least squares are autocorrelated, we used Generalized Least Squares (GLS) with seasonal AutoRegressive Integrated Moving Average (ARIMA) errors. 
	Both models obtained are GLS with ARIMA $(1, 0, 0)\times(1, 0, 0)_{11}$ errors.
	Residuals of both models satisfy all assumptions based on sample (partial) autocorrelation function, residual plot, Ljung-Box test, and Kolmogorov-Smirnov test.
	We remark that even though our model requires annual estimates, we prefer to first work with monthly data and turn those into annual data.
	This way, the predictions become more accurate, and they not only include seasonal weather patterns but also climate variability.
	
	From 12 years  of historical data, we have $m=\{1,\ldots, 252\}$ monthly data including the dependent variable $\mathrm{GPCD}_{m}$ and regressors $\mathrm{Temperature}_{m}$, $\mathrm{Precipitation}_{m}$, $\mathrm{Year}_{m}$, and binary indicator variables for each month, 
$I_1, \ldots, I_{12}$. 
For example, $\{I_{1}, \ldots, I_{12}\}=\{1, 0, \ldots,0 \}$ represents January. %
The higher-GPCD model restricts the $\mathrm{Year}_{m}$ variable to stabilize water effciciency.
First, we define residual at time ${m}$ as
$r_{m}= \mathrm{GPCD}_{m} - \textbf{X}_{m}'\vec{\bm{\beta}},$
where $\textbf{X}_{m} = [\mathrm{Temperature}_{m}, \mathrm{Precipitation}_{m}, \mathrm{Year}_{m}, I_1, \ldots, I_{12}]$ and $\vec{\bm{\beta}}=\{{\beta}_1,\ldots
, {\beta}_{15}\}$ are corresponding estimated parameters.
Let $\psi_1$ and $\Psi_1$ be estimated parameters---for autoregressive AR(1) and seasonal autoregressive SAR(1), respectively---on  the  residual. 
This time-series model with  lag operator $B$ is
$(1-\psi_1B)(1-\Psi_1B^{11})r_{m}= \varepsilon_{m},$
where $B^{t}r_{m}=r_{m-t}$ and the random noise $\varepsilon_{m}$ follows a normal distribution with mean zero and constant variance (per usual assumptions on errors). 
As a result, to predict GPCD, we use\vspace*{-0.1in}
\begin{small}
\begin{eqnarray*}\label{eq: time2}
&&\mathrm{GPCD}_{m} -\textbf{X}_{m}'\vec{\bm{\beta}}  \\
&=&\psi_1\left(\mathrm{GPCD}_{{m}-1} - \textbf{X}_{{m}-1}'\vec{\bm{\beta}}\right) + \Psi_1\left(\mathrm{GPCD}_{{m}-11} - \textbf{X}_{{m}-11}'\vec{\bm{\beta}}\right) \nonumber 
 - \psi_1\Psi_1\left(\mathrm{GPCD}_{{m}-12} - \textbf{X}_{{m}-12}'\vec{\bm{\beta}}\right).  
\end{eqnarray*}
\end{small}
The estimated parameters $\vec{\bm{\beta}}, \psi_1, \Psi_1$ are listed in Table O.2 of the Online Supplement.

	\medskip
	\noindent
	\textbf{Estimating Future Demands.}
	By the above analysis, we now have two functions to estimate average GPCD in a future month.
	These functions take as input future temperature and precipitation predictions of climate models with a given greenhouse concentration pathway.
	The predicted average GPCDs in future months are turned into average GPCDs in future years by simply considering the number of days in a month and year.
	Finally, the annual demands in a zone are estimated by multiplying the GPCDs with the predicted population of a zone. 
	We now summarize the climate models and greenhouse concentration paths used in the study.

	Bias-Corrected and Spatially Downscaled (BCSD) data from Coupled Model Intercomparison Project: Phase 5 (CMIP5) was obtained from \cite{Brekke2013}. %
	We picked the climate models listed in Table \ref{tb:climate_models} to have a good representation without overly increasing the problem size. 
	Additional climate models can be easily added to the study.%

	\begin{table}[ht]
	\small
		\centering
		\begin{tabular}{p{.7\textwidth}|r}
			\hline
			\centering{Institution} & Model \\
			\hline
			Commonwealth Scientific and Industrial Research Organization (CSIRO) and Bureau of Meteorology (BOM), Australia & CSIRO-mk-3-6-0 \\
			\hline
			\multirow{2}{*}{Geophysical Fluid Dynamics Laboratory} & GFDL-CM3 \\
			& GFDL-ESM2M \\
			\hline
			Met Office Hadley Centre & HadGEM2-ES \\
			\hline
			\multirow{3}{*}{\vbox{Atmosphere and Ocean Research Institute (The University of Tokyo), National Institute for Environmental Studies and Japan Agency for Marine-Earth Science and Technology}}  & MIROC5 \\
			& MIROC-ESM \\
			& \\
			\hline
		\end{tabular}
			\caption{
			A list of climate models used in the analysis.
		}
		\label{tb:climate_models}
	\end{table}

	Each climate model works with a given path for future greenhouse gas concentrations, called the Representative Concentration Pathway (RCP).
	Our analysis includes the four paths RCP2.6, RCP4.5, RCP6.0 and RCP8.5 adopted by IPCC \citep{pachauri2014climate}.
	RCP2.6 is an optimistic case, where concentrations are drastically reduced by mid-century.
	The paths RCP4.5 and RCP6.0 show stabilization of concentrations before and after 2060, respectively.
	Finally, RCP8.5 is the case where concentrations continue to grow quickly throughout the remainder of the century.

\begin{figure}[ht]
 	\centering
	\noindent\includegraphics[scale=0.35]{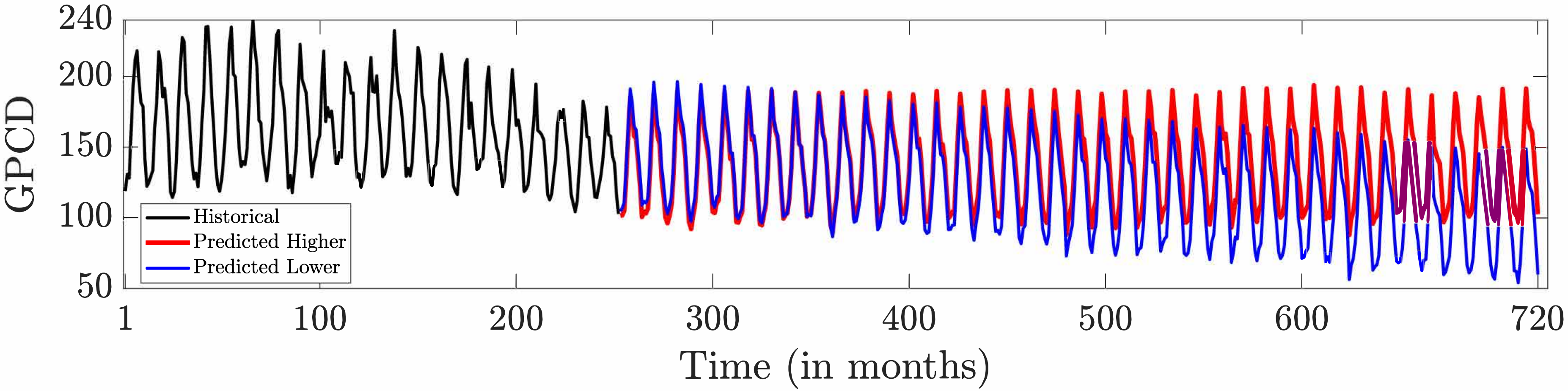}\vspace{-0.5cm}
	\caption{The lower- and higher-GPCD demand projection for the climate model CSIRO-mk-3-6-0, with greenhouse concentration pathway RCP8.5.}
	\label{fig:highlow}
\end{figure}

		The results of these predictions for one climate model, CSIRO-mk-3-6-0, with one greenhouse gas concentration pathway, RCP8.5, are shown in Figure \ref{fig:highlow}.  

	\vspace*{-0.19in}
	\subsubsection{Water Supply Prediction}
	\label{ssec: supply_prediction}

	Annual water supply to the study area is calculated by ``Allocation to Tucson $\times$ $\frac{\text{Study Area Pop.}}{\text{Tucson Pop.}}$.''
	Population estimates were discussed earlier in Section \ref{ssec: population}. 
	Below, we explain how we estimate future Colorado River water allocation to Tucson. 
	
 The conditions described in the Colorado River Compact 2007 Interim Guidelines \citep{usbor2007} dictate the Colorado River water allocation to Tucson.
	Under normal condition, Tucson Water has an annual water allocation of 144,000 af.  
	According to the compact, there are three drought conditions: Tiers 1, 2, and 3.    
	Tier 1 drought is declared if Lake Mead elevation is between 1,050--1,075 feet by end of December in a given year. %
	If so, allocation is reduced by 11.43\%. %
	Tier 2 water shortage happens when Lake Mead elevation belongs in the range $[1,025,\ 1,050)$. 
	Then, the water allocation is reduced by 14.29\%. %
	Finally, under extreme water shortage of Tier 3 (Lake Mead elevation below 1,025 feet), only 119,318 af is allocated to Tucson---a 17.14\% reduction.

	To predict the future water allocations, we used Lake Mead elevation simulations of the \cite{usbor_12,nowak}.
	We estimated the {\it nominal} probability of each condition---normal, tiers 1, 2, and 3---as the fraction of all end-of-December Lake Mead elevation simulations that satisfy a specific condition at least once during a given stage. 
	Table \ref{tb:shortage_pro} summarizes the results. 
	These simulations indicate that the chance of normal condition decreases and the chance of extreme shortage increases over the years.
	\begin{table}[ht]
	\small
		\centering
		\begin{tabular}{c|c|cccc} \hline
			\multirow{2}{*}{Stage} &  \multirow{2}{*}{Years} &  \multicolumn{4}{c}{Probability of Conditions}                                                                                                                                                      \\
			&                                                &          \begin{tabular}[c]{@{}c@{}}Normal\\ (144,000 af)\end{tabular}               & \begin{tabular}[c]{@{}c@{}}Tier 1\\ (127,541 af)\end{tabular} & \begin{tabular}[c]{@{}c@{}}Tier 2\\ (123,422 af)\end{tabular} & \begin{tabular}[c]{@{}c@{}}Tier 3\\ (119,318 af)\end{tabular} \\ \hline
			1      & \multicolumn{1}{l|}{2018}       & 1.0000  &   0.0000 & 0.0000  & 0.0000 \\
			2     & 2019--2026      & 0.6038 & 0.0817 & 0.0725 & 0.2420 \\
			3      & 2027--2034      & 0.4699 & 0.1014 & 0.0800 & 0.3488   \\
			4     & 2035--2042      & 0.3990 & 0.0854 & 0.0686 & 0.4470   \\
			5      & 2043--2050      & 0.3663 & 0.0805 & 0.0532 & 0.5000    \\   \hline      
		\end{tabular}
		\caption{
Estimated nominal probabilities of  water allotment conditions.
		}\vspace*{-0.03in}
		\label{tb:shortage_pro}
	\end{table}
	\subsubsection{Scenarios and Infrastructure Configurations}
	\label{ssec: infrastructure}
	Putting this all together, we consider the following uncertain elements at each stage:
    48 climate-related per-capita demand scenarios (=4 greenhouse gas concentrations $\times$ 6 climate models $\times$ 2 per-capita water use models, (higher-GPCD, lower-GPCD)); 2 population projections (high- and low-populations); and 4 water allotment scenarios (normal, and tier 1,2,3 droughts). 
	
	Each scenario path follows the same climate and greenhouse gas concentration pathway. This results in  384 ($=48 \times 2 \times 4)$ second-stage scenarios. Once the climate models are fixed, the other stages consider changes in population and Colorado River water allotment, resulting in 8 scenarios per node. These yield a total of 196,608 future scenarios for our study.
	We change the nominal probability of each scenario only according to its water allotment. 
	All other uncertainties are assumed to be equally likely because we do not have a preference for  climate models, population models, etc. 
	For example, a realization $\xi_2$ with normal water allotment in the second stage has the conditional probability  $q(\xi_2|\xi_{[1]})=\left(4 \times 0.6038\right)/384 =   0.0063.$

	In addition to the scenarios outlined, three infrastructure  options are considered	in the study area:  \texttt{NI} (no additional infrastructure is constructed); \texttt{WWTP} (a satellite wastewater treatment plant is constructed, capable of treating wastewater up to a nonpotable quality, for satisfying demands in its own zone and and higher zones); and \texttt{IPR} (in addition to the WWTP, an indirect potable reuse facility is constructed, which further treats water from the WWTP up to potable quality).
	Figure \ref{fig:tucson_treatment} illustrates the additional WWTP and IPR constructed in Zone C.\vspace*{-0.1in}
		
	\section{Numerical Results}   \vspace{-0.1in}
	\label{sec: result}

	Algorithm \ref{alg:decomp} is implemented in Python3.6 using the linear programming solver CPLEX 12.8.0., where we used a stopping tolerance of $\epsilon=1$E-3.
	All experiments were run on a PC with Intel Xeon Siver 4112 CPU, 2.6GHz, and 128 GB memory.
	We used python \textit{time.process\_time()} module to measure the CPU time.
	For each $\phi$-divergence, we considered the size of ambiguity set corresponding to the asymptotic confidence regions of $90\%$, $95\%$, and $99\%$ at each stage by using the $\chi^2$ value discussed in Section \ref{ssec: int_phi_divergences}. The second stage has the smallest conditional ambiguity sets and all other stages have the same, larger size\footnote{As an example, with KL divergence at 95\% confidence, $\rho_1=0.5594$ and $\rho_t=0.8792$ for  $t>1$.}.

	\subsection{Performance of the Decomposition Algorithm\note[id=JP]{Need more time to fill out the table}}
	\label{ssec: cqp_vs_decomp}
	We first computationally test the performance of Algorithm \ref{alg:decomp} with \texttt{NI} water allocation model using Hellinger distance at 95\% confidence level.
	For comparison purposes, we consider a three-stage MDRO with balanced scenario trees.
	We contrast Algorithm \ref{alg:decomp} to the CPLEX solution of SOCP formulation presented in the Online Supplement.

	Table \ref{tb:decomp_cqp} summarizes the results. 
	Column $n$ lists the total number of scenarios. 
	Remaining columns denote the running time of Algorithm \ref{alg:decomp} with single-cut, multi-cut, and SOCP in minutes, respectively. %
	Solution times reported do not include problem construction times.
	However, we note that 
	SOCP
	takes a very long time to construct with a large number of scenarios.
	
	Table \ref{tb:decomp_cqp} reveals that the decomposition algorithm is much faster than the SOCP, except for the smallest instances, and the multi-cut variant consistently outperforms the  single-cut variant.
	With 10,000 scenarios, the proposed decomposition algorithm is already more than 40 times faster than  direct solution of SOCP.
	As expected, decomposition is critical to be able to solve the water allocation problem, which is about 20 times larger than the largest problem listed in Table \ref{tb:decomp_cqp} (in terms of scenarios). It is not even possible to construct the extensive SOCP formulation of the water allocation problem without decomposition. 
	Because Algorithm \ref{alg:decomp} with multi-cut is faster than the single-cut variant, we use the multi-cut variant throughout the rest of the numerical results to analyze the water allocation problem.
	
	\begin{table}[ht]
	\small
		\centering
  \small{   \begin{tabular}{r|rr|rrr|rr|r}
\cmidrule{1-4}\cmidrule{6-9}    \multicolumn{1}{c|}{\multirow{2}[2]{*}{$n$}} & \multicolumn{2}{c|}{Algorithm 1} & \multicolumn{1}{c}{\multirow{2}[2]{*}{SOCP}} &       & \multicolumn{1}{c|}{\multirow{2}[2]{*}{$n$}} & \multicolumn{2}{c|}{Algorithm 1} & \multicolumn{1}{c}{\multirow{2}[2]{*}{SOCP}} \\
          & \multicolumn{1}{c}{single-cut} & \multicolumn{1}{c|}{multi-cut} &       &       &       & \multicolumn{1}{c}{single-cut} & \multicolumn{1}{c|}{multi-cut} &  \\
\cmidrule{1-4}\cmidrule{6-9}    4     & 0.05  & 0.05  & 0.01  &       & 2025  & 24.33 & 11.89 & 30.52 \\
    49    & 0.46  & 0.24  & 0.23  &       & 4096  & 45.28 & 29.28 & 237.35 \\
    529   & 11.33 & 4.19  & 4.22  &       & 6561  & 68.50  & 44.65 & 717.55 \\
    1,024 & 12.19 & 6.63  & 9.47  &       & 10,000 & 95.83 & 69.40  & 2976.31 \\
\cmidrule{1-4}\cmidrule{6-9}    \end{tabular}%
    }		\caption{Running time (minutes) of decomposition method (single/multi-cut) and SOCP.}\vspace*{-0.03in}
		\label{tb:decomp_cqp}
	\end{table}
   \vspace{-0.05in}

	 \vspace{-0.15in}
	\subsection{Optimal Costs and Worst-Case Probabilities by $\phi$-Divergence}
	\label{ssec: result_cost}
	We begin our analysis by comparing the optimal costs and worst-case probabilities by $\phi$-divergence. 
	
	Table \ref{tb:total_cost_any_infrastructure} lists the optimal expected costs---including operating and water-shortage costs---by infrastructure type, $\phi$-divergence, and confidence level.
	Across different confidence levels and infrastructure options, 
	the modified $\chi^2$ distance generates the highest costs followed by KL divergence and Hellinger distance.
	Burg entropy produces the lowest cost.
	The major difference between these $\phi$-divergences %
	is in the scenarios they suppressed.
	Recall that a scenario is suppressed if its optimal worst-case probability is zero.
	The modified $\chi^2$ distance suppressed scenarios for every confidence level tested, and it consistently suppressed the lower-GPCD  scenarios.
	Especially,  with 99\% confidence it suppressed \textit{all} the lower-GPCD  scenarios, including both high- and low-population cases.
	In other words, scenarios with relatively low demands are ignored.
	KL divergence and Hellinger distance maintain an ``all-or-nothing'' approach to suppressing scenarios, but confidence levels of 90--99\% are not high enough to induce the suppressing behavior for this problem for any stage at any node. Instead, they put low optimal worst-case probabilities on the lower-GPCD  scenarios.
	
	For this specific problem, we find that if a scenario $\omega_{t+1}$ is suppressed ($p^{\omega_{t+1}|\omega_t*}=0$) with modified
    $\chi^2$ distance, we have an order:  $p^{\omega_{t+1}|\omega_t*}$ with KL is less than that with  Hellinger, which is less then that of Burg entropy for the same scenario $\omega_{t+1}$.
	We conjecture that above worst-case probability order on low-cost scenarios explains the cost order in Table \ref{tb:total_cost_any_infrastructure}. 
	In additional tests, we observed that other instances of this problem and different problems in other domains do not show this behavior. 
	\begin{table}[ht]
	\small
		\centering
		\begin{tabular}{l|rrr|rrr|rrr}
    \midrule
          & \multicolumn{3}{c|}{\texttt{NI}} & \multicolumn{3}{c|}{\texttt{WWTP}} & \multicolumn{3}{c}{\texttt{IPR}} \\
    \midrule
    \multicolumn{1}{c|}{$\phi$-divergence} & \multicolumn{1}{c}{90\%} & \multicolumn{1}{c}{95\%} & \multicolumn{1}{c|}{99\%} & \multicolumn{1}{c}{90\%} & \multicolumn{1}{c}{95\%} & \multicolumn{1}{c|}{99\%} & \multicolumn{1}{c}{90\%} & \multicolumn{1}{c}{95\%} & \multicolumn{1}{c}{99\%} \\
    \midrule
     Modified  $\chi^2$   & 458.67 & 461.10 & 465.26 & 433.45 & 435.80 & 439.84 & 400.59 & 402.75 & 406.44 \\
     Kullback-Leibler     & 453.54 & 457.74 & 463.01 & 429.05 & 433.13 & 438.29 & 395.81 & 399.56 & 404.24 \\
     Hellinger            & 444.89 & 449.82 & 457.17 & 421.34 & 425.98 & 432.97 & 387.71 & 392.20 & 398.86 \\
     Burg                 & 437.00 & 441.53 & 448.74 & 414.28 & 418.54 & 425.34 & 380.32 & 384.45 & 391.02 \\
    \bottomrule
    \end{tabular}%
    \caption{
			Optimal expected costs  (in \$ million) over 2018--2050 for each infrastructure configuration.
		}
		\label{tb:total_cost_any_infrastructure}
	\end{table}
    
    In the rest of this section we highlight our results using the KL divergence at 95\% confidence level. Other $\phi$-divergences and confidence levels produce similar results.

	\subsection{Comparison of Climate Models and Greenhouse Gas Concentration Paths}
	\label{ssec: result_climate}
	
	Table \ref{tb:pworst_model_emission_kl} presents the total optimal worst-case probabilities assigned to each climate model and greenhouse gas concentration path by MDRO.
	Among the concentration paths, 
	MDRO assigns the highest probability to the highest concentration path RCP8.5 and lowest probability to the lowest concentration path RCP2.6.
	And among the climate models, it assigns the highest probability to GFDL-CM3 followed by HadGEM2-ES, which tend to generate higher temperatures than other climate models (see Table O.3 in the Online Supplement).
	These results indicate that MDRO tends to put higher probabilities to models that together generate higher temperatures. Higher temperatures lead to increased water demands, and thus increased costs. 
	This way, MDRO induces a risk-averse behavior, protecting against the more frequent water shortages associated with these scenarios.\vspace*{-0.1in}

	\begin{table}[ht]
	\small
		\centering
		\begin{tabular}{c|l|cccc|c}
&	& \multicolumn{4}{c|}{ Greenhouse Gas Concentration Paths} & \\
& & RCP2.6 & RCP4.5 & RCP6.0 & RCP8.5 & (all) \\ \hline
\parbox[t]{2mm}{\multirow{6}{*}{\rotatebox[origin=c]{90}{Climate Models}}}
          & CSIRO &  0.0399 & 0.0436 & 0.0375 &0.0408 &  0.1618 \\
          & GFDL-CM3 &  0.0417 & 0.0440 & 0.0430 & 0.0453& 0.1740  \\
          & GFDL-ESM2M & 0.0387 & 0.0384 & 0.0394 &0.0394 & 0.1559 \\
          & HadGEM2-ES & 0.0422 & 0.0421 & 0.0427 &0.0455 & 0.1725 \\
          & MIROC5 & 0.0393 & 0.0413 & 0.0401 &0.0431& 0.1638 \\
          & MIROC-ESM-CHEM & 0.0410 & 0.0428 & 0.0434 & 0.0450 & 0.1722\\ \hline
          & (all) &  0.2428 & 0.2522 & 0.2461 & 0.2591 & 1 \\
\end{tabular}%
\caption{
			Optimal probabilities for each climate model and concentration path (KL, 95\%).
		}
		\label{tb:pworst_model_emission_kl}
	\end{table}
	\vspace*{-0.1in}

	\subsection{Evaluation of Additional Decentralized Infrastructures}
	\label{ssec: options}
	We now use our MDRO model to help governing agencies in their infrastructure decisions.%
	
\begin{figure}[ht]
		\centering
			\includegraphics[scale=0.19]{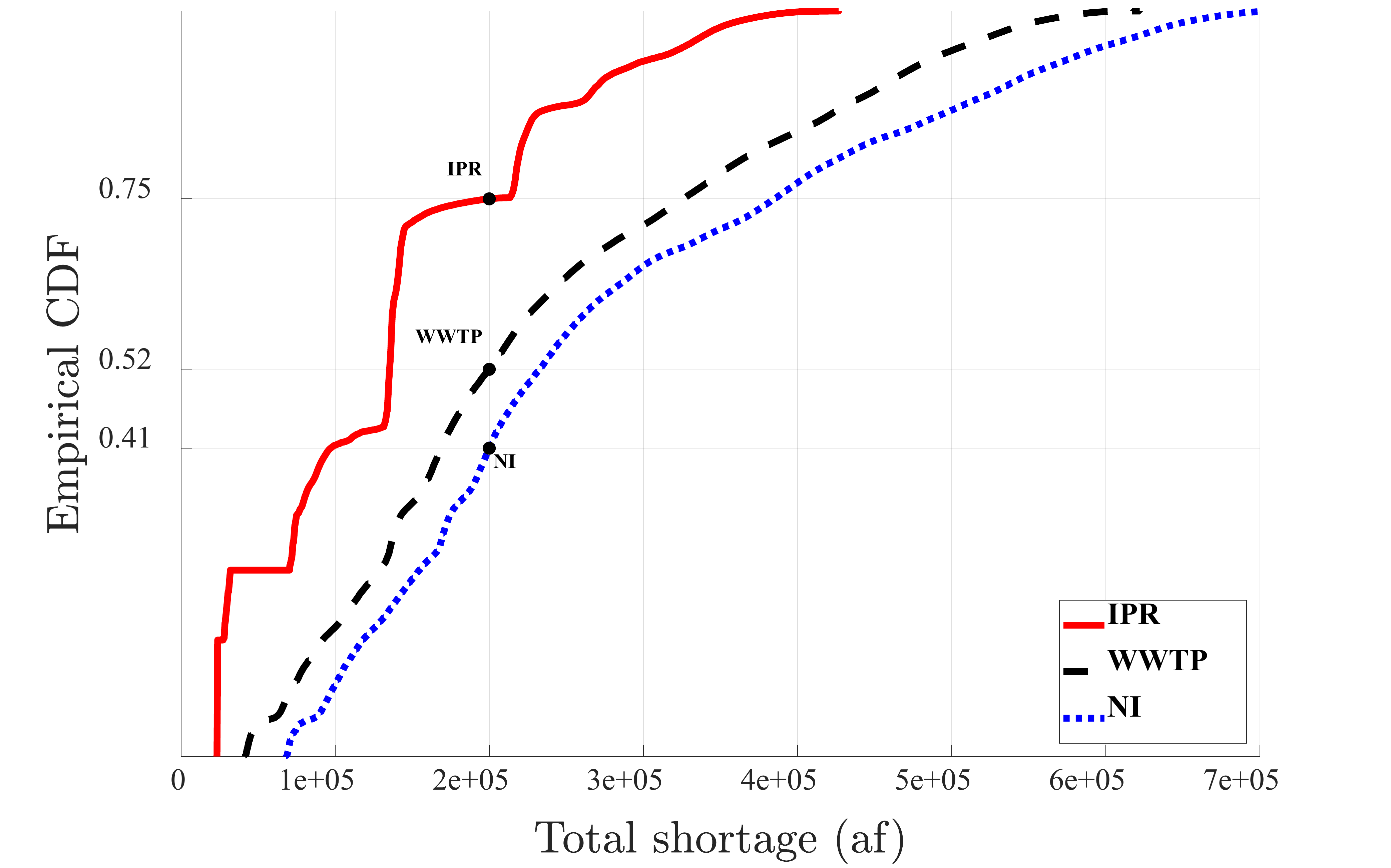}\vspace*{-0.2in}
			\caption{
				Empirical CDF of total shortage over the 33-year study period for each infrastructure configuration  (KL, 95\%).
			}\vspace*{-0.05in}
			\label{fig:cdf_shortage_combined}
	\end{figure}

	\noindent {\textbf{Water Shortage.}}
	One of the most important advantages of decentralized water treatment is that it increases water reuse. 
	With this additional water supply, water shortage is decreased. 
	Therefore, we examine the effect of additional infrastructure on water shortage first. 
	Figure \ref{fig:cdf_shortage_combined} depicts the empirical Cumulative Distribution Function (CDF) of the total shortage (in af) for each infrastructure configuration using the nominal distribution.
	The CDFs of \texttt{WWTP} and \texttt{IPR} are always above that of \texttt{NI}. 
	This means that they are preferable to \texttt{NI} regarding shortage, and they stochastically dominate not having decentralized infrastructures in the area.
	For example, the nominal probability that total shortage is less than or equal to $200,000$ af is 0.41, 0.52 and 0.75 for \texttt{NI}, \texttt{WWTP} and \texttt{IPR}, respectively. 
	Looking at the highest CDFs, \texttt{IPR} provides the most substantial reduction in water shortages, followed by 
	\texttt{WWTP} as a somewhat distant second. 
	
\begin{figure}[!th]
	\centering
	\includegraphics[width=0.8\textwidth]{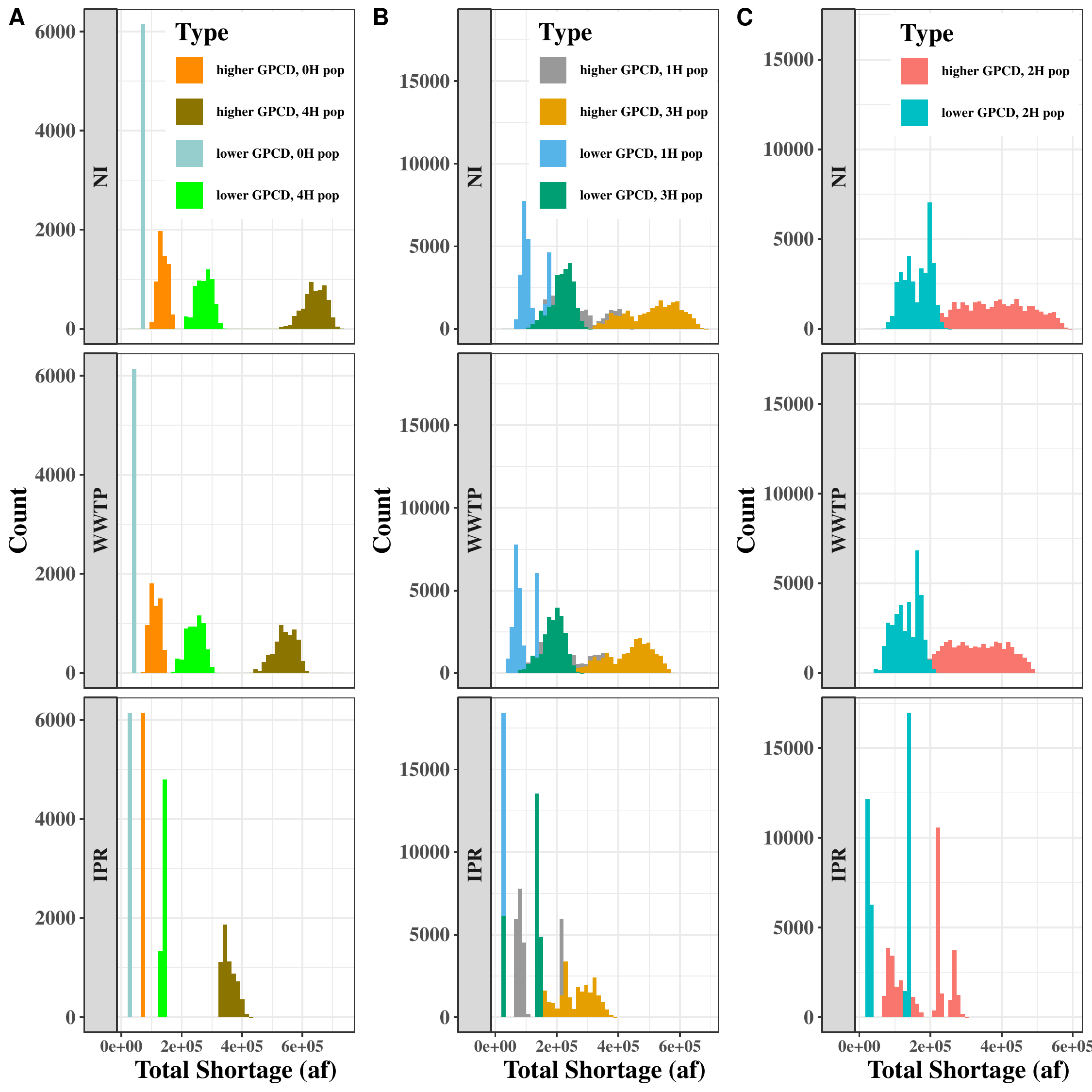}\vspace*{-0.07in}
			\caption{
				Histogram of total shortage amount over the 33-year study period for each GPCD, population and infrastructure configuration (KL, 95\%).
			}\vspace*{-0.03in}
			\label{fig:shortage_frequency_combined}
	\end{figure}

	Figure \ref{fig:shortage_frequency_combined} depicts the total shortage amount over the 33-year study period, broken down by  infrastructure,  %
	GPCD demand, and population categories. \textit{``i"H} indicates \textit{``i"} stages have high-population scenarios. For example, \textit{0H pop} indicates all stages have low-population scenarios and \textit{4H pop} means all stochastic stages have high-population scenarios (first stage is deterministic).
	A decentralized \texttt{WWTP}
	provides a considerable reduction in shortage severity---especially in higher-GPCD scenarios with at least one high-population stage.  
	The \texttt{IPR} facility substantially decreases shortages of all categories. Especially, the \texttt{IPR} is effective in reducing the \textit{extreme} shortages. 
	
	\vspace*{-0.01in}

	The result of this analysis clearly shows the value of \texttt{IPR}---and to a lesser extent the value of \texttt{WWTP}---in reducing water shortages. 
	We examine their economic value next.
	\medskip

	\noindent {\textbf{Cost-Benefit and Break-Even Analysis.}}
	Table \ref{tb:total_cost_any_infrastructure} reveals that \texttt{WWTP} consistently decreases the operating cost by about \$24.61 million, and the \texttt{IPR} facility reduces the cost by an additional \$23.57  (or total \$58.18) million over the 33-year time span.
	This is mainly due to the reduced water shortages. 
	An earlier analysis by Tuscon Water indicates that the \texttt{WWTP} and the \texttt{IPR} facility, if constructed, would cost \$55 and \$119(=55+64) million, respectively.
	As a result, the additional facilities will {\it not} pay for themselves over the planning period. %

	So far we used a shortage cost of \$800/af. 
	This led us to examine break-even shortage costs that balance the construction costs with operational savings.  
	We find that, compared to \texttt{NI}, 
	\texttt{WWTP} breaks even at \$1,765/af and the additional \texttt{IPR} facility breaks even at \$1,256/af shortage cost.
	The break-even cost for \texttt{IPR} is lower because it drastically lowers shortages.
	Assuming a satellite \texttt{WWTP} in the area is already built, \texttt{IPR}  breaks-even at \$1,100/af.
	These results imply that the increased operation cost of \texttt{IPR} plus its higher construction cost is far lower than the benefits it provides when the shortage costs are increased.\vspace*{-0.09in}

	\subsection{Discussion}\vspace*{-0.09in}
	The above analysis reveals that per-capita water demand (measured in GPCD) is the main driver of water shortages among the categorized uncertainties (Figure \ref{fig:shortage_frequency_combined}).
	The high shortages in Figure \ref{fig:shortage_frequency_combined} occur at higher-GPCD scenarios when at least one stage has a high-population scenario.
	The implication of this result is twofold. 
	First, water conservation efforts and technologies could have a drastic effect in the area. 
	Second, next to water-shortage costs,  a decision on  building additional infrastructures largely depends on GPCD.
	It should not be overlooked that GPCD in itself depends on climate models and greenhouse concentration paths.
	Therefore, the final construction decision needs to consider the impact of climate scenarios.
	
	Another important conclusion of this study is that \texttt{IPR} is the best option in terms of shortages (total shortage amounts and break-even shortage costs) especially as shortage cost increases; see, e.g., Figure \ref{fig:cdf_shortage_combined}. However, public opinion and long-term health effects should be considered before constructing an IPR facility because  drinking heavily treated wastewater  has not yet been supported by the public \citep{Kerri2013Drinking, Laura}.\vspace*{-0.1in}
	
	\section{Concluding Remarks}\vspace*{-0.01in}
	\label{sec: concl_water}
 In real life the true distribution governing the random parameters is never fully known. This issue becomes more serious for multistage problems. This concern motivated us to consider the so-called distributionally robust approach. In particular, we built MDRO models with conditional ambiguity sets of distributions on a given scenario tree by staying sufficiently close to nominal conditional distributions using $\phi$-divergences. We devised a nested Benders decomposition algorithm to solve this class of problems. The algorithm provides a significant generalization of the subgradient-based decomposition that was earlier used for the simpler mean-CVaR case. 
 Next, we put the MDRO modeling and solution techniques to use by solving a real-life water allocation problem under the uncertainties of climate, population, and Colorado River availability, among others. And we evaluated the value of decentralized infrastructures. 
 To the best of our knowledge this is the first study of MDRO for managing water resources under climate uncertainty.  

The results of this paper can lay the foundation for studying new algorithms and models. For instance, studying a sampling-based version, e.g., stochastic dual dynamic programming \citep{pereira_pinto_91}, of the decomposition method discussed in the paper would be valuable to approximately solve larger models. 
As our study is strategic in nature, the need for large detail is low. However, for many operational models, it would be valuable to have larger-size models that can only be approximately solved via sampling.
Many nontrivial research questions arise in this case including  how to generate samples within the algorithm to speed convergence and  how to stop the algorithm with rigorous stopping criteria. These merit further, rigorous investigations. %

In terms of water allocation model, incorporation of water quality and especially health impacts of \texttt{IPR} would be valuable, but this requires long-term studies. Modeling the dependence of water price to climate events also merits further study.  
Lastly, in addition to the water allocation problem discussed in this paper, the investigated MDRO modeling and solution techniques have the potential to make an impact on other problems with ambiguous  time-dynamic uncertainties, e.g., that arise in other environmental, energy management, and financial problems.

\vspace*{-0.1in}
\begin{small}
\baselineskip0.15in
\bibliographystyle{apalike}
\bibliography{ref}

\end{small}

\end{document}